\documentclass[final]{siamart1116}

\usepackage{amsmath,amssymb,euscript,latexsym,graphicx}
\usepackage{bbm,color,amstext,wasysym,subfig,parskip,balance}
\graphicspath{{./},{./figures/}}

\newtheorem{thm}{Theorem}

\newtheorem{remark}[thm]{Remark}

\newcommand{\be}{\begin{equation}}
\newcommand{\ee}{\end{equation}}
\newcommand{\mR}{{\mathbb R}}

\newcommand{\cD}{{\mathcal D}}
\newcommand{\cE}{{\mathcal E}}
\newcommand{\cF}{{\mathcal F}}
\newcommand{\cG}{{\mathcal G}}

\newcommand{\cV}{{\mathcal V}}

\newcommand{\cW}{{\mathcal W}}
\newcommand{\diag}{\operatorname{diag}}

\newcommand{\la}{\langle}
\newcommand{\ra}{\rangle}

\definecolor{grey}{rgb}{0.6,0.6,0.6}
\definecolor{lightgray}{rgb}{0.97,.99,0.99}

\begin{document}
\title{Vector-Valued Optimal Mass Transport}
\author{Yongxin Chen\thanks{Y.\ Chen is with the Department of Medical Physics, Memorial Sloan Kettering Cancer Center, NY; email: chen2468@umn.edu}\and
Tryphon T. Georgiou\thanks{T.\ T. Georgiou is with the Department of Mechanical and Aerospace Engineering, University of California, Irvine, CA; email: tryphon@uci.edu}\and Allen Tannenbaum\thanks{A.\ Tannenbaum is with the Departments of Computer Science and Applied Mathematics \& Statistics, Stony Brook University, NY; email: allen.tannenbaum@stonybrook.edu}}

\maketitle

\begin{abstract}
We introduce the problem of transporting vector-valued distributions. In this, a salient feature is that mass may flow between vectorial entries as well as across space (discrete or continuous).
The theory relies on a first step taken to define an appropriate notion of optimal transport on a graph. The corresponding distance between distributions is readily computable via convex optimization and provides a suitable generalization of Wasserstein-type metrics.
Building on this, we define Wasserstein-type metrics on vector-valued distributions supported on continuous spaces as well as graphs. Motivation for developing vector-valued mass transport is provided by applications such as color image processing, multi-modality imaging, polarimetric radar, as well as network problems where resources may be vectorial.
\end{abstract}

\section{Introduction}
The theory of Monge-Kantorovich optimal mass transport theory has witnessed a fast pace of new developments; see \cite{RacRus98,Vil03} for extensive lists of references. These contributions were driven by a multitude of applications in physics, geosciences, economics, and probability. Some of the notable advances include the concept of displacement interpolation \cite{Mcc97}, links to the geometry of spaces \cite{BakEme85, LotVil09, Ott01, StuRen05},  and a fluid dynamic reformulation \cite{BenBre00}. In our own work, image analysis and spectral analysis of time series provided starting points (e.g., \cite{GKT09,HakTanAng04}) and, more recently, problems in stochastic control, quantum information, and matrix-valued distributions \cite{CheGeoPav14e,JiaNinGeo12,NinGeoTan15,TGT}. The present paper continues the work of \cite{CheGanGeoTan17,CheGeoTan16} by proposing a transportation problem for vector-valued distributions.

A salient feature of vector-valued distributions is the possibility of the transfer of ``mass'' from one vectorial entry to another. Physical examples include color image scenes where the vectorial distribution captures color intensities, which may continuously shift with lighting conditions. Alternatively, polarimetric data provides an analogous example where mass represents power at different polarizations detected at the locations of a sensor array. As another example, the flow of mass between vector entries may represent mutation of coexisting population species.

The proposed framework may have far-reaching consequences in, for example, combining genomic and proteomic networks, and in general, fusion of vectorial data supported on a graph. While in some of these examples the total mass may not be preserved, in the present work we will restrict our attention to the case where it is. Thus, we seek suitable continuity equation that allows trading off mass between vectorial entries of a distribution, on a continuous or a discrete space (graph), and develop a geometric framework that would allow constructing geodesic flows between snapshots of such distributions.

In order to formulate transport between vectorial entries, we begin with a new notion of transport on weighed undirected graphs in the spirit of Erbar and Maas \cite{ErbMaa12}. A starting point in \cite{ErbMaa12} is to devise a suitable continuity equation for probability measures on the nodes of a weighted graph (Markov chain). The formulation in our paper differs from that in \cite{ErbMaa12}, and the corresponding transport problem has the advantage of being reducible to one of convex optimization. Both \cite{ErbMaa12} and our formulation were inspired by the Benamou-Brenier theory \cite{BenBre00}, where the OMT with quadratic cost is recast as the problem to minimize flow kinetic energy (i.e., an action integral). The present work builds on \cite{CheGanGeoTan17,CheGeoTan16} extending the Wasserstein theory to densities and mass distributions on more general spaces. Having as a first step a Benamou-Brenier theory on graphs, the methodology allows us to define a notion of vector-valued transport and corresponding distance between vector-valued densities on discrete or continuous spaces. As with the (weighted) graph case, the transport distance that we define on vector-valued densities may be reduced to a convex optimization problem.

We now outline the remainder of this note. In Section~\ref{sec:pre}, we sketch needed background from the classical theory of optimal mass transport that motivates our generalization. In Section~\ref{sec:graph}, we describe the proposed Wasserstein-2 metric on an undirected weighted graph. Further, we remark on Wasserstein-1 type of metric on a weighted graph. In Sections~\ref{sec:vector} and \ref{sec:vectorgraph}, we formulate the new Wasserstein distance on vector-valued densities that are supported, first on the Euclidean spaces and then on graphs. In Section~\ref{sec:examples}, we give several examples illustrating the idea of a vector-valued optimal mass transport, and finally we conclude in Section~\ref{sec:conclusion} with an outline of possible applications of the theory and future research directions.

\section{Preliminaries on optimal mass transport}\label{sec:pre}

The mass transport problem was first formulated by Gaspar Monge in
1781, and concerned finding the optimal way, in the sense of minimal
transportation cost, of moving a pile of soil from one site to
another. This problem was given a modern formulation in the work of
Kantorovich in the form of a linear program and it is now known as the {\em Monge--Kantorovich
problem.} See \cite{EvaGan99, Kan48, RacRus98, Vil03} for all details as well as extensive lists of references.

Herein, we focus mainly the case where the transportation cost is quadratic in the distance. The respective optimization problem
	\begin{equation}\label{eq:scalarprimal}
		W_p (\rho_0,\rho_1)^p:=\inf_{\pi\in \Pi(\rho_0,\rho_1)} \int_{\mR^N\times\mR^N} \|x-y\|^p \pi(dx,dy)
	\end{equation}
for $p=2$, where $\|\cdot\|$ denotes Euclidean distance and $\Pi(\rho_0,\rho_1)$ represents the set of all couplings
between two $\rho_0$ and $\rho_1$ non-negative probability density functions on $\mR^N$ (i.e., the set of joint probability distributions having $\rho_0$ and $\rho_1$ as respective marginals), defines the so-called Wasserstein-$2$ distance between the two densities, or more generally, between measures.

In this case, where the cost is quadratic (i.e., $p=2$), the transport problem admits a dynamic reformulation \cite{BenBre00} that is especially powerful and the space of densities $\cD :=\{\rho \ge 0\mid\int \rho\, dx =1 \}$ admits, essentially, a Riemannian structure \cite{Ott01}. The Benamou-Brenier reformulation identifies the Wasserstein-2 distance with the integral of the kinetic energy (action integral) along a geodesic flow that links the two marginals, namely,
\begin{subequations}
\begin{equation}\label{eq:minimization}
W_2 (\rho_0,\rho_1)^2= \inf_{\rho,v}
\int \int_0^1 \rho(t,x) \|v(t,x)\|^2 \, dt \, dx
\end{equation}
over all time varying densities $\rho$ and vector fields $v$ satisfying the continuity equation and boundary conditions
\begin{eqnarray}
&&\frac{\partial \rho}{\partial t} + \nabla_x \cdot (\rho v) = 0,
\label{eqcont}\\
&&\rho(0, \cdot) = \rho_0, \; \rho(1, \cdot) = \rho_1 \label{eq:marginals}.
\end{eqnarray}
\end{subequations}
Interestingly, when expressed in terms of density $\rho$ and flux $m=\rho v$, the minimization problem in \eqref{eq:minimization} becomes convex while (\ref{eqcont}-\ref{eq:marginals}) turn into linear constraints.
For the optimal pair $(\rho,v)$, the vector field turns out to be the gradient $v=\nabla_x g$ of a function $g$, hence, ``rot-free''.
Vector fields $v$ of this form can be identified with tangent directions of $\cD$, i.e., elements of the tangent space
\[
T_\rho \cD \cong \{\delta: \int \delta =0 \},
\]
as follows. Under suitable assumptions on
differentiability for $\rho \in \cD$, and $\delta \in T_\rho \cD$, we solve the Poisson equation
\begin{equation} \label{poisson}
\delta= -\nabla_x \cdot (\rho \nabla_x g)
\end{equation}
to obtain a convex function $g_\delta$ and thereby the vector field $\nabla g_\delta$.
In this way the space $\cD$ can be endowed with a Riemannian structure (see \cite{Ott01,Vil03}) via
\begin{equation} \label{inner} \la \delta_1, \delta_2 \ra_\rho := \int \rho \nabla_x g_{\delta_1} \cdot \nabla_x g_{\delta_2}dx \end{equation}
which has the aforementioned kinetic energy interpretation. This inner product induces precisely the Wasserstein distance as geodesic distance between the two marginals in \eqref{eq:marginals}.

\begin{remark}[Wasserstein-1 metric]\label{remark:W1} Interestingly, the case of linear cost (i.e., $p=1$), can also be cast as the flux minimization problem
\begin{equation*}\label{eq:minimization2}
W_1 (\rho_0,\rho_1)= \inf_{m}
 \int  \|m\| \, dx
\end{equation*}
\begin{eqnarray*}
 \rho_1-\rho_0+ \nabla_x \cdot m &=& 0,
\label{eqcont2}
\end{eqnarray*}
as it can be shown that the flux $m=\rho v$ remains invariant with time \cite{EvaGan99}.
For this case, there is also an alternative expression through a dual formulation
\begin{equation*}\label{eq:minimizationW1second}
W_1 (\rho_0,\rho_1)= \sup_{f}\{ \int f (\rho_1-\rho_0)dx \mid \|\nabla_x f\| \leq 1\}.
\end{equation*}
in terms of test functions; see \cite{EvaGan99,RacRus98,Vil03,LiYinOsh16}.
\end{remark}

\section{Wasserstein metric on weighted graphs}\label{sec:graph}

Following the Benamou-Brenier viewpoint to Wasserstein distances our first task is to develop an analogous notion of transportation distances on graphs. To this end, we consider a connected, positively weighted, undirected graph $\cG=(\cV, \cE,\cW)$ with $n$ nodes labeled as $i$, with $1\le i \le n$, and $m$ edges.  We consider the set of probability masses on $\cG$ that we will denote by $\cD$; an element $\rho \in \cD$ may be regarded as a column vector $(\rho_1, \ldots, \rho_n)^T$, with $\rho_i \ge 0$ for $1 \le i \le n$ and
\[
\sum_{i=1}^n \rho_i =1.
\]
We denote the (open) interior of $\cD$ by $\cD_+$.

The standard heat equation on $\cG$,
	\[
		\dot\rho=L\rho=-DWD^T \rho,
		\]
	where $L, D, W=\diag \{w_1, \cdots, w_m\}$ are the graph-Laplacian, incidence, and weight  matrices, respectively, can also be written in the more familiar (from calculus)
	\begin{equation}\label{eq:heat}
			\dot\rho=\Delta_\cG \rho
	\end{equation}
by defining
\[  \Delta_\cG := -\nabla_\cG^*\nabla_\cG,
\]
where
	\[
		\nabla_\cG: \mR^n \rightarrow \mR^m, ~ x \mapsto W^{1/2}D^T x
	\]
denotes the gradient operator and
	\[
		\nabla_\cG^*:  \mR^m \rightarrow \mR^n, ~y \mapsto DW^{1/2} y
	\]
denotes its dual. More generally, if we
let the entries of $u(t) \in \mR^m$ represent flux along respective edges, we can express the continuity equation in the form
\be \label{eq:cont1}
		\dot \rho -\nabla_\cG^* u=0.
\ee
Evidently, the flux $u=-\nabla_\cG \rho$ gives \eqref{eq:heat}.
Also note that since the row vector consisting of all 1's lies in the left kernel of the incidence matrix, mass is preserved by (\ref{eq:cont1}).

To carry out our program, we need to express the flux $u$ in the form of a momentum ``$\rho v$'' as in \cite{BenBre00}. However, the flux is supported on the edges $\cE$ of the graph whereas the mass is supported on the set of nodes $\cV$, the two sets having different dimensions. In order to overcome this difficulty in a natural manner, we choose to
associate the flux along an edge with the mass at the source in the two end-points.
More specifically, the flux along an edge $e_k=(i,j)$, with source $i$ and sink $j$, consists of two parts. A part that flows out of node $i$, and another that flows in opposite direction out of node $j$. Thus, we define a flux $u_k=\rho_i v_k$ out of $i$ and another, $\bar u_k=\rho_j\bar v_k$ out of $j$, and represent the total flux as the superposition $\rho_j  v_k-\rho_i \bar v_k$, while restricting the rates $v_k, \bar v_k$ to be nonnegative. Thus, our
 \textbf{\emph{continuity equation}} for rates $v, \bar v\in \mR_+^m$ becomes
	\be \label{eq:cont2}
	\dot\rho-\nabla_\cG^*((D_2^T\rho)\circ v-(D_1^T\rho)\circ \bar v)=0,
	\ee
where $\circ$ denotes entry-wise multiplication of two vectors. The matrix $D_1$ is the portion of the incidence matrix $D$ containing $1$'s (sources), and $D_2=D_1-D$ (sinks). In other words, $D_1^T\rho$ is the mass at the source of an edge, and $D_2^T\rho$ is the mass at the sink of an edge. {\em The dependence of the flux in (\ref{eq:cont2}) on $\rho$ ensures that the entries of $\rho$ remain positive while the fact that the kernel of $D$ contains the vector with all ones ensures that the total mass is preserved as well.}

For notational convenience we use $\mu\in\cD_+$ and $\nu\in\cD_+$ (instead of $\rho_0$ and $\rho_1$) to denote the starting and ending mass on nodes. We now define the transport distance between $\mu$ and $\nu$ as follows:
	\begin{eqnarray}\label{eq:firstform}
	W_{2,a}(\mu,\nu)^2 :=&& \inf_{\rho,v,\bar v} \int_0^1 \left\{v^T ((D_2^T\rho)\circ v)+
	\bar v^T((D_1^T\rho)\circ \bar v)\right\}dt
				\\&& \dot\rho-\nabla_\cG^*((D_2^T\rho)\circ v-(D_1^T\rho)\circ \bar v)=0,\nonumber
				\\&& v\ge 0, ~~\bar v\ge 0,\nonumber
				\\&& \rho(0)=\mu,~~\rho(1)=\nu.\nonumber
	\end{eqnarray}
It is easy to see that at each time instant, for each $k$, at most one of the $v_k$ and $\bar v_k$ is nonzero.
In a similar manner as in the Benamou-Brenier program, \eqref{eq:firstform} can be recast in the form of a convex optimization problem in (momentum) variables $u=(D_2^T\rho)\circ v, \bar u=(D_1^T\rho)\circ \bar v$,
	\begin{eqnarray}\label{eq:graphconv1}
	W_{2,a}(\mu,\nu)^2 =&& \inf_{\rho,
	u,\bar u} \int_0^1 \left\{u^T \diag(D_2^T\rho)^{-1}u+
	\bar u^T \diag(D_1^T\rho)^{-1}\bar u\right\}dt
				\\&& \dot\rho-\nabla_\cG^*(u- \bar u)=0,\nonumber
				\\&& u\ge 0, ~~\bar u\ge 0, \label{eq:graphconv3}\nonumber
				\\&& \rho(0)=\mu,~~\rho(1)=\nu.\nonumber
	\end{eqnarray}

It is straightforward to see that the right hand side in \eqref{eq:graphconv1} is in general positive and vanishes only when $\mu=\nu$. It is also straightforward to see that $W_{2,a}$ satisfies the triangle inequality. However, in general, $W_{2,a}(\mu,\nu)\neq W_{2,a}(\nu,\mu)$, therefore $W_{2,a}$ is only a quasimetric. Yet, it endows $\cD_+$ with a Finsler metric type structure
	\begin{eqnarray}
		W_{2,a}(\mu,\mu+\delta)^2=&&\inf_{\rho,v,\bar v} v^T ((D_2^T\mu)\circ v)+
	\bar v^T((D_1^T\mu)\circ \bar v
		\\&&  \delta-\nabla_\cG^*((D_2^T\mu)\circ v-(D_1^T\mu)\circ \bar v)=0,\nonumber
		\\&& v\ge 0, ~~\bar v\ge 0,\nonumber
	\end{eqnarray}
for small perturbation $\delta$, and in this, $(\cD_+, W_{2,a})$ becomes a length space. In fact, $W_{2,a}$ has a very nice ``geodesic'' property. Indeed, if $\rho(\cdot)$ is the mass distribution as a function of time obtained by solving  \eqref{eq:graphconv1}, then
	\begin{equation}
		W_{2,a}(\rho(s),\rho(t))=(t-s)W_{2,a}(\mu,\nu)
	\end{equation}
for any $0\le s<t\le 1$. Finally, $W_{2,a}$ can be extended to $\cD$, the closure of $\cD_+$, by continuity.

Naturally, one can symmetrize $W_{2,a}$ in the obvious way, by adopting as our metric
$$
\max\{W_{2,a}(\mu,\nu), W_{2,a}(\nu,\mu)\},$$ which can then be computed by solving two convex optimization problems. (A similar remark holds for the metrics $W_{2,b}$ and $W_{2,c}$ defined later on.)

An alternative way to symmetrize $W_{2,a}$ is to replace the cost function \eqref{eq:graphconv1} with
	\[
		\int_0^1 \frac12\left\{u^T (\diag(D_2^T\rho)^{-1}+\diag(D_1^T\rho)^{-1})u
		+\bar u^T (\diag(D_2^T\rho)^{-1}+\diag(D_1^T\rho)^{-1})\bar u\right\}dt.
	\]
Since the cost terms for $u$ and $\bar u$ are symmetric, we can combine the two and drop the nonnegativity requirement  $u\ge 0, \bar u\ge 0$, to obtain
	\begin{eqnarray}\label{eq:sgraphconv}
	\hat W_{2,a}(\mu,\nu)^2 :=&& \inf_{\rho\in\cD_+,u} \int_0^1\left\{u^T (\diag(D_2^T\rho)^{-1}
	+\diag(D_1^T\rho)^{-1})u\right\}dt
				\\&& \dot\rho-\nabla_\cG^* u=0,\nonumber
				\\&& \rho(0)=\mu,~~\rho(1)=\nu.\nonumber
	\end{eqnarray}
Positive entries of $u$ represent flow from sources to sinks, while negative entries, flow from sinks to sources. This (symmetric) metric induces a Riemannian type structure on $\cD_+$, akin to that of standard optimal transport theory on Euclidean spaces \cite{Vil03}.

\begin{remark}We point out that similar notions was recently considered in \cite{ErbMaa12,ChoHuaLiZho12,SolRusGuiBut16,ChoLiZho17}. We arrived at the above formulation independently and from a different starting point.
\end{remark}

\begin{remark}[Gradient flow of entropy] The gradient flow of the entropy functional on probability mass distribution on a graph $\cG$ with respect to $\hat W_{2,a}$ is given by
\begin{equation}\label{eq:nonlinearheat}
\dot{\rho} = -\nabla_\cG ^* (A(\rho)^{-1} \nabla_\cG \log \rho)
\end{equation}
where $A(\rho):= \diag(D_2^T\rho)^{-1} +\diag(D_1^T\rho)^{-1}$. It represents a nonlinear heat-like equation, to be contrasted with the linear heat equation derived in \cite{ErbMaa12}. To see \eqref{eq:nonlinearheat}, compute the derivative of the  entropy functional
$
S(\rho)= -\sum_{i=1}^n \rho_i \log \rho_i
$
along a curve $\rho(t), t \in [0,1]$ in $\cD_+$,
\begin{eqnarray*}
-\dot{S}(\rho) &=& \sum_{i=1}^n \dot{\rho}_i\log \rho_i + \sum_{i=1}^n \dot{\rho_i}\\
&=& \langle \dot{\rho}, \log \rho \rangle \;\;\;\mbox{ \rm (since for each $t,$ $\sum_{i=1}^n {\rho_i}(t) =1$)}\\
&=& \langle \nabla_\cG ^*u, \log \rho \rangle  \;\;\;\mbox{ \rm (since $\dot\rho=\nabla_\cG^* u$)}\\
&=&\langle u, \nabla_\cG \log \rho \rangle \end{eqnarray*}
and observe that the direction of steepest ascent is along
$u = -A(\rho)^{-1} \nabla_\cG \log \rho,$ from which \eqref{eq:nonlinearheat} follows.
\end{remark}

\begin{remark}[Wasserstein-1 distance on graphs]\label{sec:w1}
Following up on Remark, \ref{remark:W1} we sketch a Benamou-Brenier type reformulation of $W_1$-distances on graphs. Assuming that the entries of $c=(c_1,\cdots,c_m)^T$ are edge weights representing the cost of moving a unit mass across, a $W_1$ distance between mass distributions $\mu,\nu$ can be defined as the solution of the min-cost flow problem
	\begin{eqnarray}\label{eq:W1cost}
	W_1(\mu,\nu) =&& \min_{u} c^T |u|
				\\&& \nu-\mu-Du=0.\nonumber
	\end{eqnarray}
Alternative, in terms of ``test vectors'' $f$,
	\begin{equation*}\label{eq:W1}
		W_1(\mu,\nu)= \max_{f} \left\{f^T(\nu-\mu) ~\mid~ \|\nabla_\cG f \|_\infty \le 1\right\},
	\end{equation*}
having the dual
	\begin{equation*}
		W_1(\mu,\nu)=\min_{\hat{u}} \left\{\|\hat u\|_1 ~\mid~ \nu-\mu-\nabla_\cG^* \hat{u} =0\right\}.
	\end{equation*}
This coincides with \eqref{eq:W1cost} by taking $c_i=1/\sqrt{w_i}$, $1\le i\le m$, and $u= W^{1/2}\hat u$.
Finally, we point out that the above has an action minimization formulation following \cite{Leo16}:
	\begin{subequations}
	\begin{eqnarray*}
	W_1(\mu,\nu) =&& \inf_{
	\rho,
	v,\bar v} \int_0^1 \left\{c^T((D_2^T\rho)\circ v)+
	c^T ((D_1^T\rho)\circ\bar v)\right\}dt
				\\&& \dot\rho-D((D_2^T\rho)\circ v-(D_1^T\rho)\circ \bar v)=0,
				\\&& v\ge 0, ~~\bar v\ge 0,
				\\&& \rho(0)=\mu,~~\rho(1)=\nu,
	\end{eqnarray*}
	\end{subequations}
which assumes the following convex recast
	\begin{subequations}
	\begin{eqnarray*}
	W_1(\mu,\nu) =&& \inf_{
	\rho,u,\bar u} \int_0^1 \left\{c^Tu+
	c^T \bar u\right\}dt
				\\&& \dot\rho-D(u-\bar u)=0,
				\\&& u\ge 0, ~~\bar u\ge 0,
				\\&& \rho(0)=\mu,~~\rho(1)=\nu.
	\end{eqnarray*}
	\end{subequations}
\end{remark}

\section{Vector-valued densities \& transport} \label{sec:vector}

We now turn to the main theme of our paper, the introduction of a Wasserstein type metric between vector-valued densities. A vector-valued density
$\rho=[\rho_1,\rho_2,\cdots,\rho_M]^T$ on $\mR^N$, or on a discrete space, may represent a physical entity that can mutate or be transported between alternative manifestations, e.g., power reflected off a surface at different frequencies or polarizations. While the total power may be invariant (under some lighting conditions), the proportion of power at different frequencies or polarization may smoothly vary with viewing angle. As another example consider the case where the entries of $\rho$ represent densities of different species, or particles, and allow for the possibility that mass transfers from one species to another, i.e., between entries of $\rho$. Thus, in general, we postulate that transport of vector-valued quantities captures flow across space as well as between entries of the density vector. We introduce an OMT-inspired geometry that allows us to express a continuity and quantify transport cost for such vectorial distributions.

We begin by considering a vector-valued density $\rho$ on $\mR^N$, i.e., a map from $\mR^N$ to $\mR^M_+$ such that
	\[
		\sum_{i=1}^M\int_{\mR^N} \rho_i(x)dx=1.
	\]
To avoid proliferation of symbols we denote the set of all vector-valued densities and its interior again by $\cD$ and $\cD_+$, respectively. We refer to the entries of $\rho$ as representing density or mass of species/particles that can mutate between one another while maintaining total mass.
The dynamics are captured by the following \textbf{\emph{continuity equation}}:
	\begin{equation}\label{eq:continuityvector}
		\frac{\partial\rho_i}{\partial t}+\nabla_x\cdot(\rho_i v_i)-\sum_{j\neq i} (\rho_j w_{ji}-\rho_i w_{ij})=0,
		~~\forall i=1,\ldots,M.
	\end{equation}
Here $v_i$ is the velocity field of particles $i$ and $w_{ij}\ge 0$ is the transfer rate from $i$ to $j$.
Equation \eqref{eq:continuityvector} allows for the possibility to mutate between each pair of entries. More generally, mass transfer may only be permissible between specific types of particles and can be modeled by a graph $\cF=(\cV_1,\cE_1,\cW_1)$. Thus,  \eqref{eq:continuityvector} corresponds to the case where $\cF$ is a complete graph with all weights equal to $1$. For general $\cF$ the continuity equation is
	\begin{equation}\label{eq:continuitygeneral}
		\frac{\partial \rho}{\partial t} +\nabla_x\cdot(\rho\circ v)-
		\nabla_\cF^* ( (D_2^T \rho)\circ w-(D_1^T\rho)\circ \bar w) = 0.
	\end{equation}
Note here $\rho$ denote the vector $[\rho_1,\rho_2,\cdots,\rho_M]^T$ and likewise for $v, w, \bar w$.

Given $\mu,\nu\in\cD_+$, we formulate the optimal mass transport:
	\begin{eqnarray}\label{eq:omt}
	W_{2,b}(\mu,\nu)^2&& := \inf_{
	\rho,v,w,\bar w}
	\int_0^1\int_{\mR^N}\left\{ v^T (\rho\circ v)+\gamma[w^T ((D_2^T\rho)\circ w)\right.\\
	&&\left.\hspace*{4.7cm}\nonumber +
	\bar w^T((D_1^T\rho)\circ \bar w)]\right\} dx dt
	\\&& \frac{\partial \rho}{\partial t} +\nabla_x\cdot(\rho\circ v)-
		\nabla_\cF^* ( (D_2^T \rho)\circ w-(D_1^T\rho)\circ \bar w) = 0\nonumber
	\\&& w(t,x) \ge 0, ~~\bar w(t,x) \ge 0, ~~\forall t,x \nonumber
	\\&& \rho(0,\cdot)=\mu(\cdot), ~ \rho(1,\cdot)=\nu(\cdot).\nonumber
	\end{eqnarray}
The coefficient $\gamma>0$ specifies the relative cost between transporting mass in space and trading mass between different types of particles. When $\gamma$ is large, the solution reduces to independent OMT problems for the different entries to the degree possible. As with $W_{2,a}$, it can be shown that $W_{2,b}$ is a quasi-metric in that it satisfies the triangle inequality and positivity, but is not symmetric. Also, $W_{2,b}$ has the geodesic property
	\begin{equation}
		W_{2,b}(\rho(s,\cdot),\rho(t,\cdot))=(t-s)W_{2,b}(\mu,\nu)
	\end{equation}
for  $0\le s< t\le 1$, assuming $\rho$ is the optimal flow for \eqref{eq:omt}.

Setting $p=\rho\circ w\ge 0, \bar p=\rho\circ \bar w\ge 0$ and $u=\rho\circ v$, we establish that \eqref{eq:omt} is equivalent to the convex optimization problem
	\begin{eqnarray}\label{eq:omtcvx1}
	&&\inf_{\rho,
	u,p,\bar p}
	\int_0^1\int_{\mR^N}\left\{ u^T\diag(\rho)^{-1} u+\gamma [p^T \diag(D_2^T\rho)^{-1}p+
	\bar p^T \diag(D_1^T\rho)^{-1}\bar p]\right\} dx dt
	\\&& \frac{\partial\rho}{\partial t}+\nabla_x\cdot u-\nabla_\cF^* (p-\bar p)=0\nonumber
	\\&& p(t,x) \ge 0, ~~\bar p(t,x) \ge 0, ~~\forall t,x \nonumber
	\\&& \rho(0,\cdot)=\mu(\cdot), ~ \rho(1,\cdot)=\nu(\cdot).\nonumber
	\end{eqnarray}
Again, as in \eqref{eq:sgraphconv}, one can define a Riemannian like metric on $\cD_+$ by symmetrizing the above, which leads to
	\begin{eqnarray}\label{eq:vectorOMTsymmetry}
	\hat W_{2,b}(\mu,\nu)^2 =&&\inf_{\rho,u,p}
	\int_0^1\int_{\mR^N}\left\{ u^T\diag(\rho)^{-1} u\right.\\\nonumber &&\left.+\gamma p^T (\diag(D_2^T\rho)^{-1}+
	\diag(D_1^T\rho)^{-1})p\right\} dx dt
	\\&& \frac{\partial\rho}{\partial t}+\nabla_x\cdot u-\nabla_\cF^* \,p=0\nonumber
	\\&& \rho(0,\cdot)=\mu(\cdot), ~ \rho(1,\cdot)=\nu(\cdot).\nonumber
	\end{eqnarray}
\begin{remark}[Wasserstein-1 distance for vector-valued densities]
With continuity equation \eqref{eq:continuitygeneral}, in the same spirit as in Remark \ref{remark:W1}, it is straightforward to define Wasserstein-1 distance for vector-valued densities as
	\begin{eqnarray*}
	W_1 (\mu,\nu)&=& \inf_{u,p}\int\{ \|u\| +\gamma \|p\|\} \, dx
 	\\&&\nu-\mu+\nabla_x\cdot u-
		\nabla_\cF^*\,p= 0,
	\end{eqnarray*}
whose dual is clearly
	\begin{eqnarray*}
	W_1 (\mu,\nu)&=& \sup_f \int\{ f^T(\nu-\mu)\}dx
	\\&& \|\nabla_x f\| \le 1,~~ \|\nabla_\cF f\| \le \gamma.
	\end{eqnarray*}
\end{remark}

\section{Vector-valued mass transport on graphs}\label{sec:vectorgraph}
We finally consider vector-valued mass transport on graphs. A vector-valued mass distribution on graph $\cG=(\cV,\cE,\cW)$ (with $n$ nodes and $m$ edges) is a $M$-tuple $\rho=(\rho_1,\cdots,\rho_M)$ with each $\rho_i=(\rho_{i,1},\cdots,\rho_{i,n})^T$ being a vector in $\mR_+^n$ such that
	\[
	\sum_{i=1}^M\sum_{k=1}^n \rho_{i,k}=1.
	\]
That is, each entry $\rho_i$, for $i\in\{1,\ldots,M\}$, is a vector with nonnegative $n$-entries representing, e.g., color intensity for the $i$-th color, at the node corresponding to the respective entry.
We denote the set of all non-negative vector-valued mass distributions with $\cD$ and its interior with $\cD_+$. Combining \eqref{eq:cont2} and \eqref{eq:continuityvector} we obtain the continuity equation
	\begin{equation}
		\dot\rho-\nabla_\cG^*((D_2^T\rho)\circ v-(D_1^T\rho)\circ \bar v)
		-\nabla_\cF^* ((D_2^T\rho)\circ w-(D_1^T\rho)\circ \bar w)=0.
	\end{equation}
The problem of transporting vector-valued mass on a graph is conceptually simpler as it reduces essentially to a scalar mass situtation. Indeed, we can view the vector-valued mass as a scalar mass distribution on $M$ identical layers of the graph $\cG$ where the same nodes at different layers are connected through a graph $\cF$. The two velocity fields $v, w$ represent mass transfer within the same layer and between different layers, respectively.

Following our earlier program, given two marginal densities $\mu, \nu\in\cD_+$, we define their Wasserstein distance as
	\begin{eqnarray*}
	W_{2,c}(\mu,\nu)^2 := &&\inf_{\rho,
	v,\bar v, w, \bar w}
	\int_0^1\left\{v^T ((D_2^T\rho)\circ v)+
	\bar v^T((D_1^T\rho)\circ \bar v)+\gamma[w^T ((D_2^T\rho)\circ w)+
	\bar w^T((D_1^T\rho)\circ \bar w)]\right\} dt
	\\&& \dot\rho-\nabla_\cG^*((D_2^T\rho)\circ v-(D_1^T\rho)\circ \bar v)
		-\nabla_\cF^*((D_2^T\rho)\circ w-(D_1^T\rho)\circ \bar w)=0,\nonumber
	\\&& w\ge 0, ~\bar w \ge 0, ~v \ge 0, ~\bar v \ge 0,
	\\&& \rho(0)=\mu, ~ \rho(1)=\nu.
	\end{eqnarray*}
In the same way as before, the above has a convex reformulation
	\begin{eqnarray*}
	&&\inf_{\rho,
	u,\bar u,p,\bar p}
	\int_0^1\left\{u^T \diag(D_2^T\rho)^{-1}u+
	\bar u^T\diag(D_1^T\rho)^{-1} \bar u+\gamma[p^T \diag(D_2^T\rho)^{-1}p+
	\bar p^T\diag(D_1^T\rho)^{-1} \bar p]\right\} dt
	\\&& \dot\rho-\nabla_\cG^*(u-\bar u)
		-\nabla_\cF^* (p-\bar p)=0,\nonumber
	\\&& p \ge 0,~\bar p \ge 0, ~u \ge 0, ~\bar u\ge 0,
	\\&& \rho(0)=\mu, ~ \rho(1)=\nu.
	\end{eqnarray*}
The same method as in \eqref{eq:sgraphconv} gives rise to a symmetric Riemannian type metric $\hat W_{2,c}(\mu,\nu)^2$ provided by the solution of
	\begin{eqnarray*}
	&&\inf_{\rho,
	u,p}
	\int_0^1\left\{ u^T (\diag(D_2^T\rho)^{-1}+
	\diag(D_1^T\rho)^{-1})u+\gamma p^T (\diag(D_2^T\rho)^{-1}+
	\diag(D_1^T\rho)^{-1})p\right\} dt
	\\&& \dot\rho-\nabla_\cG^*\,u
		-\nabla_\cF^*\, p=0,\nonumber
	\\&& \rho(0)=\mu, ~ \rho(1)=\nu. \nonumber
	\end{eqnarray*}

\section{Examples} \label{sec:examples}
In this section, we present two examples. The first one is an academic example to illustrate the idea of vector-valued optimal mass transport. In the second example, we apply our framework to color image processing problems.

\subsection{Interpolation of $1$-D densities}
\begin{figure}[h]
\centering
\subfloat[$\mu$]{\includegraphics[width=0.40\textwidth]{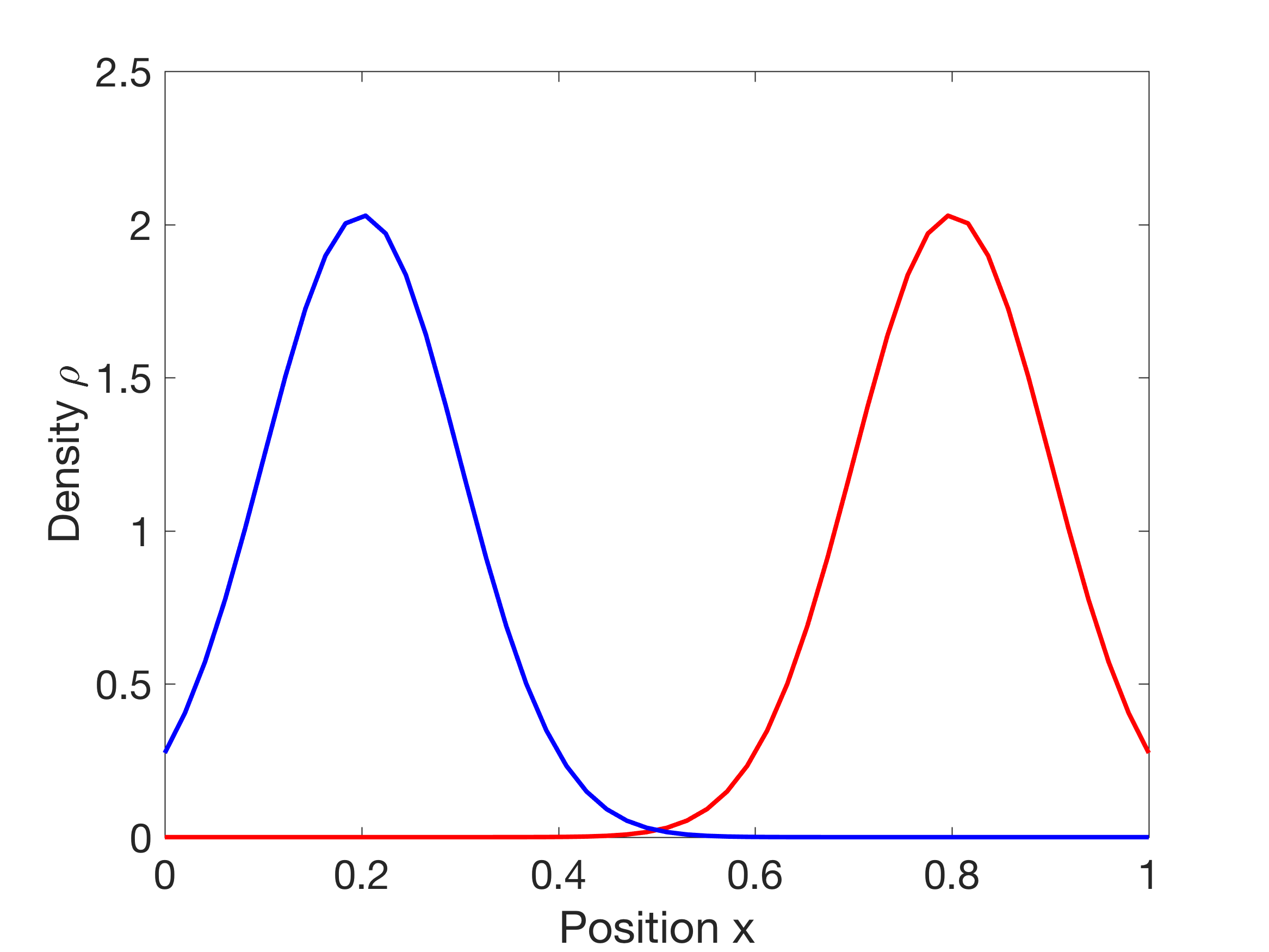}}
\subfloat[$\nu$]{\includegraphics[width=0.40\textwidth]{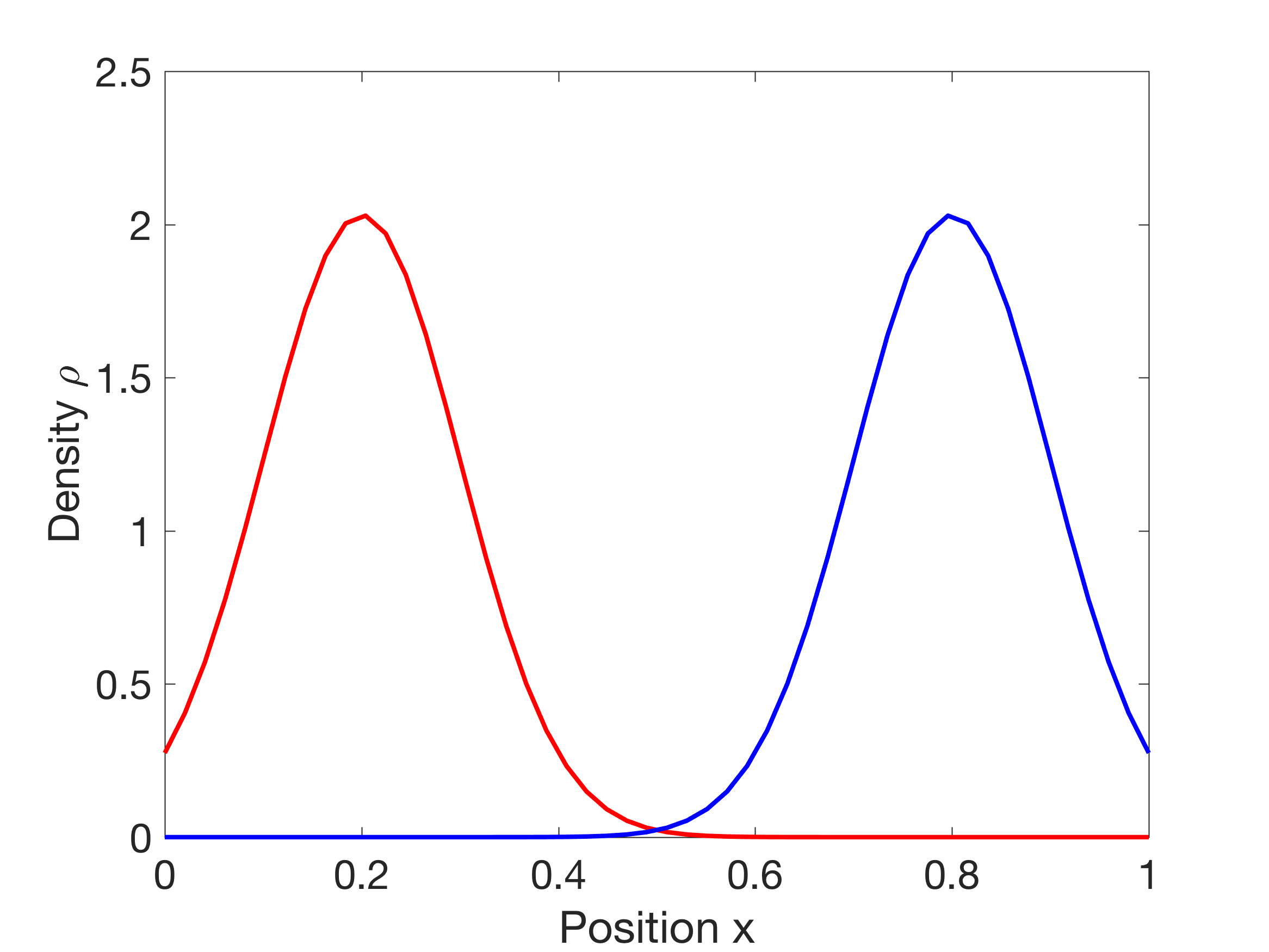}}
    \caption{Marginal distributions}
    \label{fig:marginals}
\end{figure}

We consider vector-valued densities with two components on the real line (interval $[0,\, 1]$). The two marginal densities $\mu$ and $\nu$ are displayed in Fig \ref{fig:marginals} with the two colors (red and blue) denoting the two components.

We solve the symmetric vector-valued transport problem \eqref{eq:vectorOMTsymmetry} for several different values of $\gamma$. As to the numerical implementation, we first discretize the space interval $[0,\, 1]$ to convert it into a vector-valued transport problems on graphs, which is essentially \eqref{eq:sgraphconv}. Then we discretize the time dimension with staggered grids. In particular, we discretize the time interval into $n_t$ subintervals. Then the densities take value at time points $0, 1/n_t, \cdots, (n_t-1)/n_t, 1$ while the fluxes take values at time points $1/(2n_t),3/(2n_t),\cdots,(2n_t-1)/(2n_t)$. We refer the reader to \cite{CheKaoHabGeoTan17} for more details about the convex optimization algorithm used for the examples in this paper.

\begin{figure}[h]
\centering
\subfloat[$\gamma=0.0001$]{\includegraphics[width=0.40\textwidth]{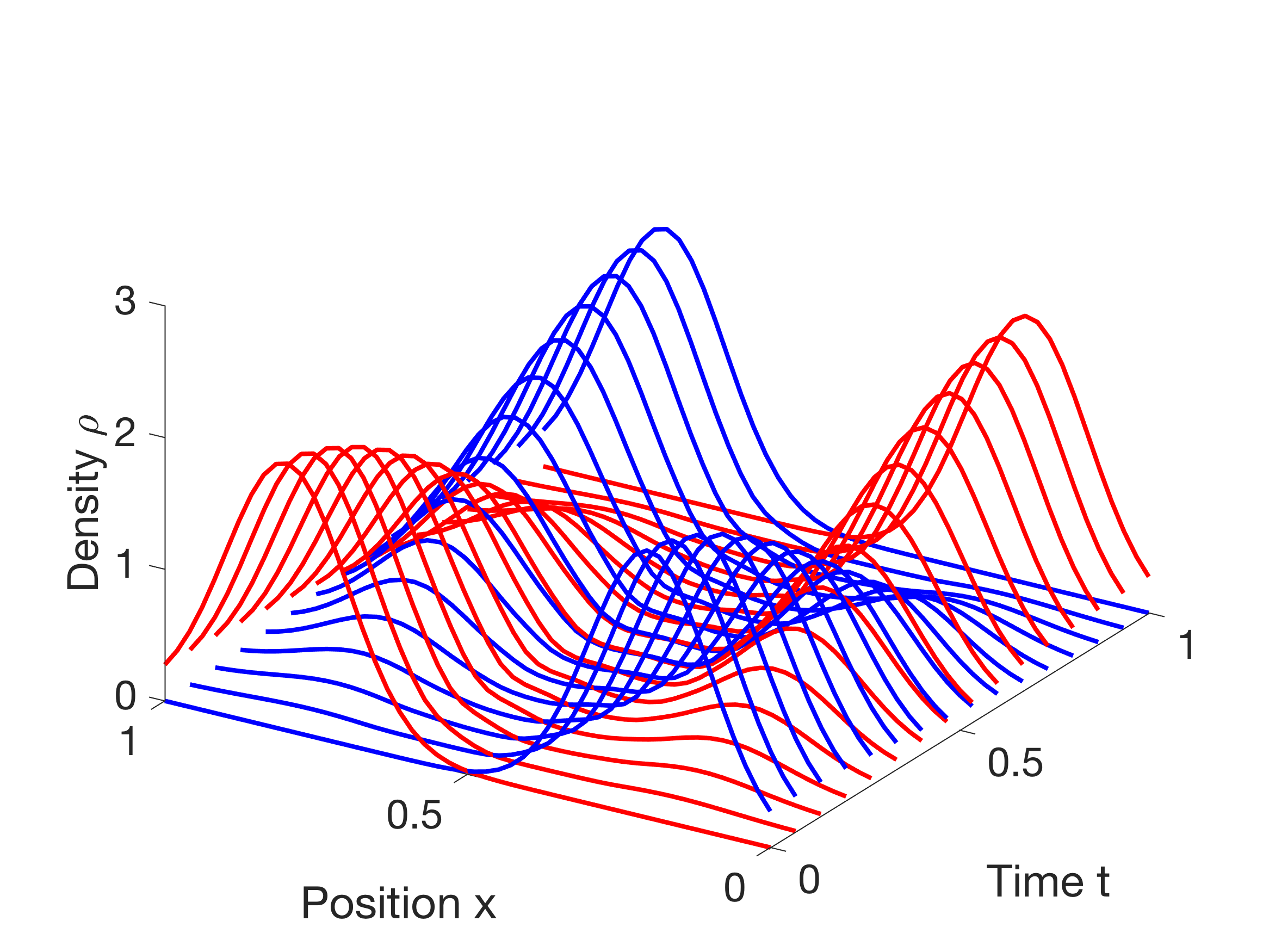}}
\subfloat[$\gamma=0.01$]{\includegraphics[width=0.40\textwidth]{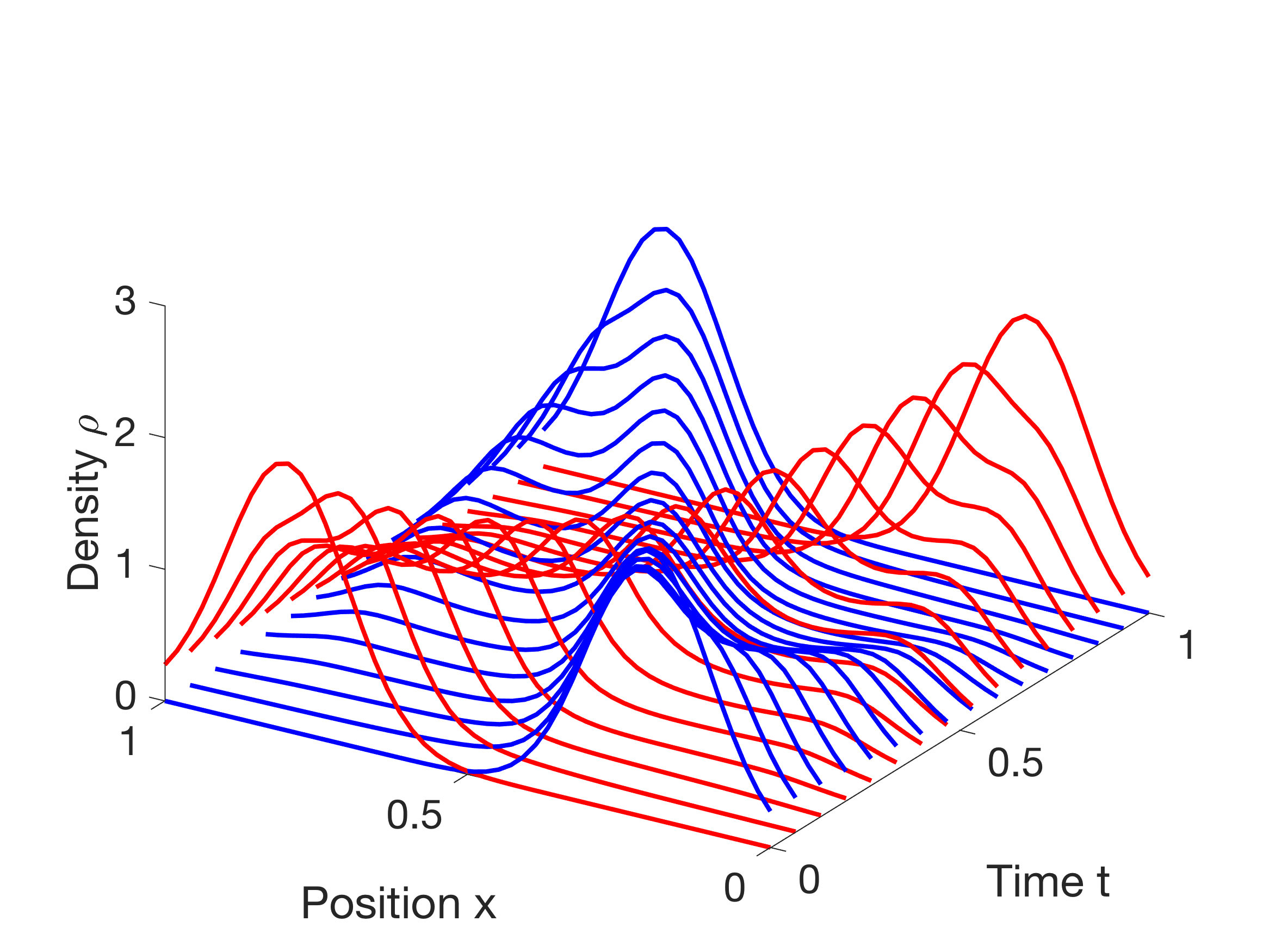}}
\\
\subfloat[$\gamma=1$]{\includegraphics[width=0.40\textwidth]{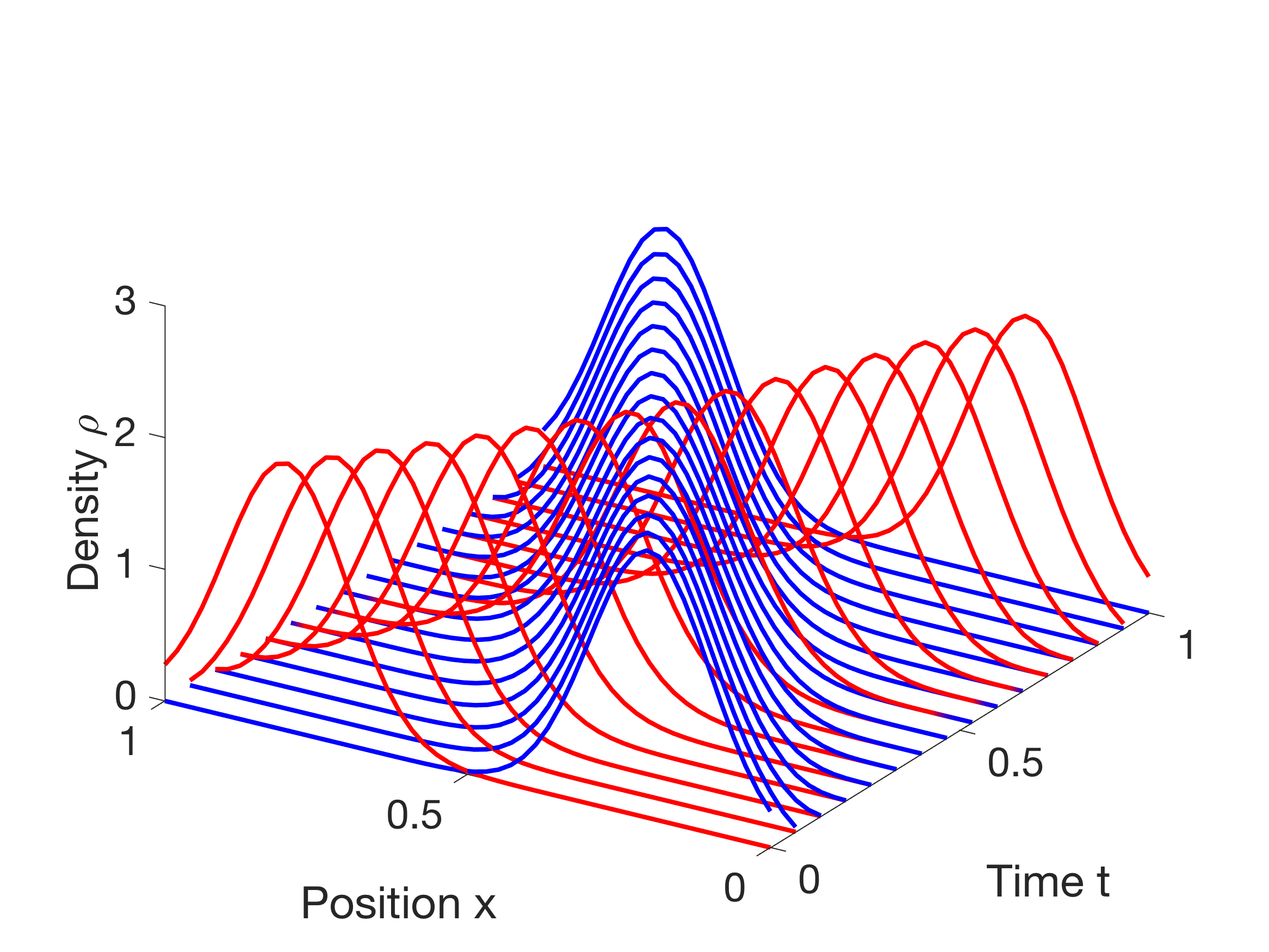}}
\subfloat[$\gamma=100$]{\includegraphics[width=0.40\textwidth]{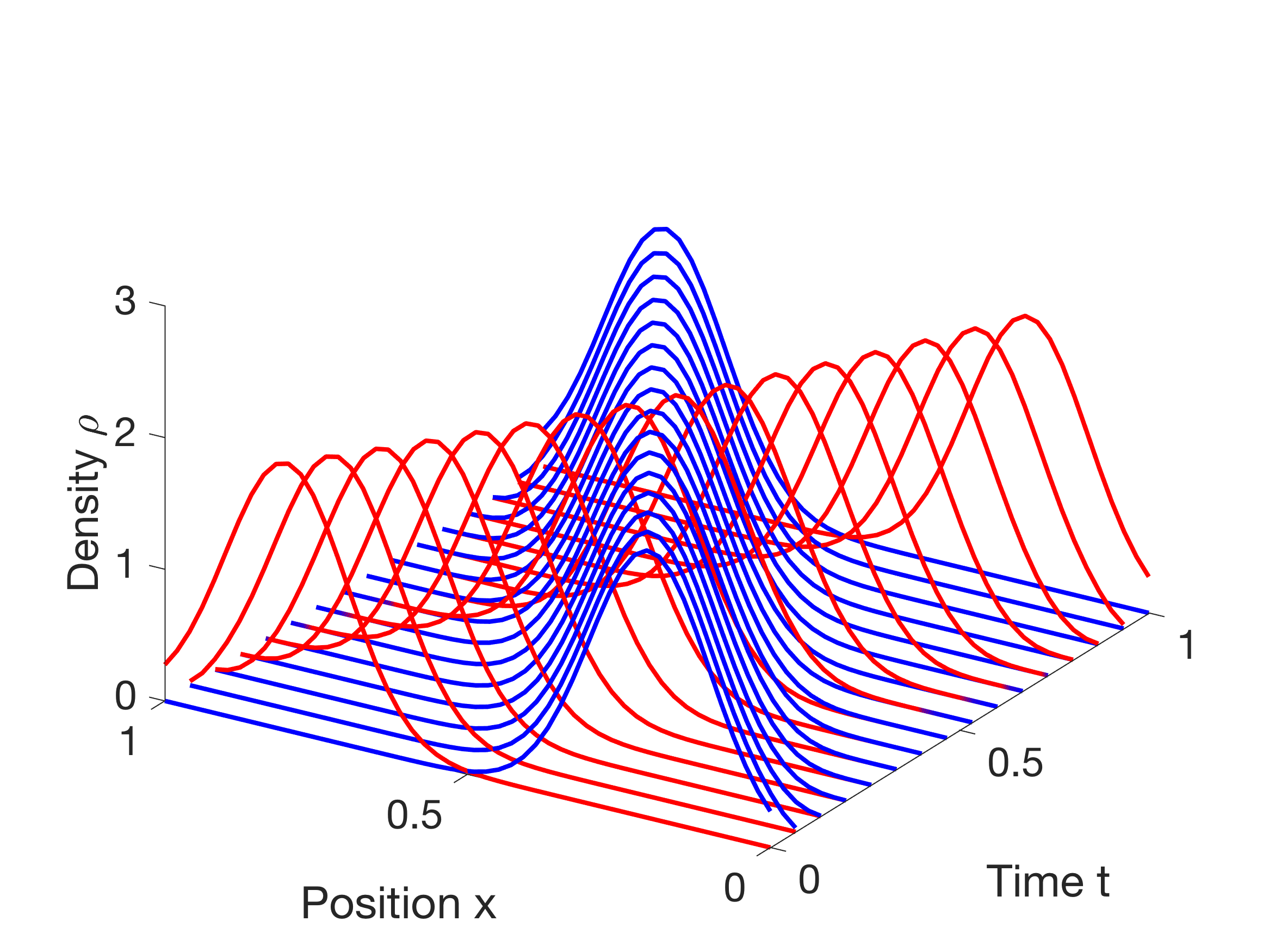}}
    \caption{Density flow with different $\gamma$ values}
    \label{fig:flow}
\end{figure}

The results are depicted in Figure~\ref{fig:flow}. As can be seen, for large $\gamma$, the solution tends to have two independent transport plan as the cost of transferring between the two different masses is high. In contrast, when $\gamma$ is small, the solution prefers transferring than transporting, since the cost of transferring between masses is low.

\subsection{Interpolation of color images}
\begin{figure}[h]
\centering
\subfloat[$\rho_0$]{\includegraphics[width=0.40\textwidth]{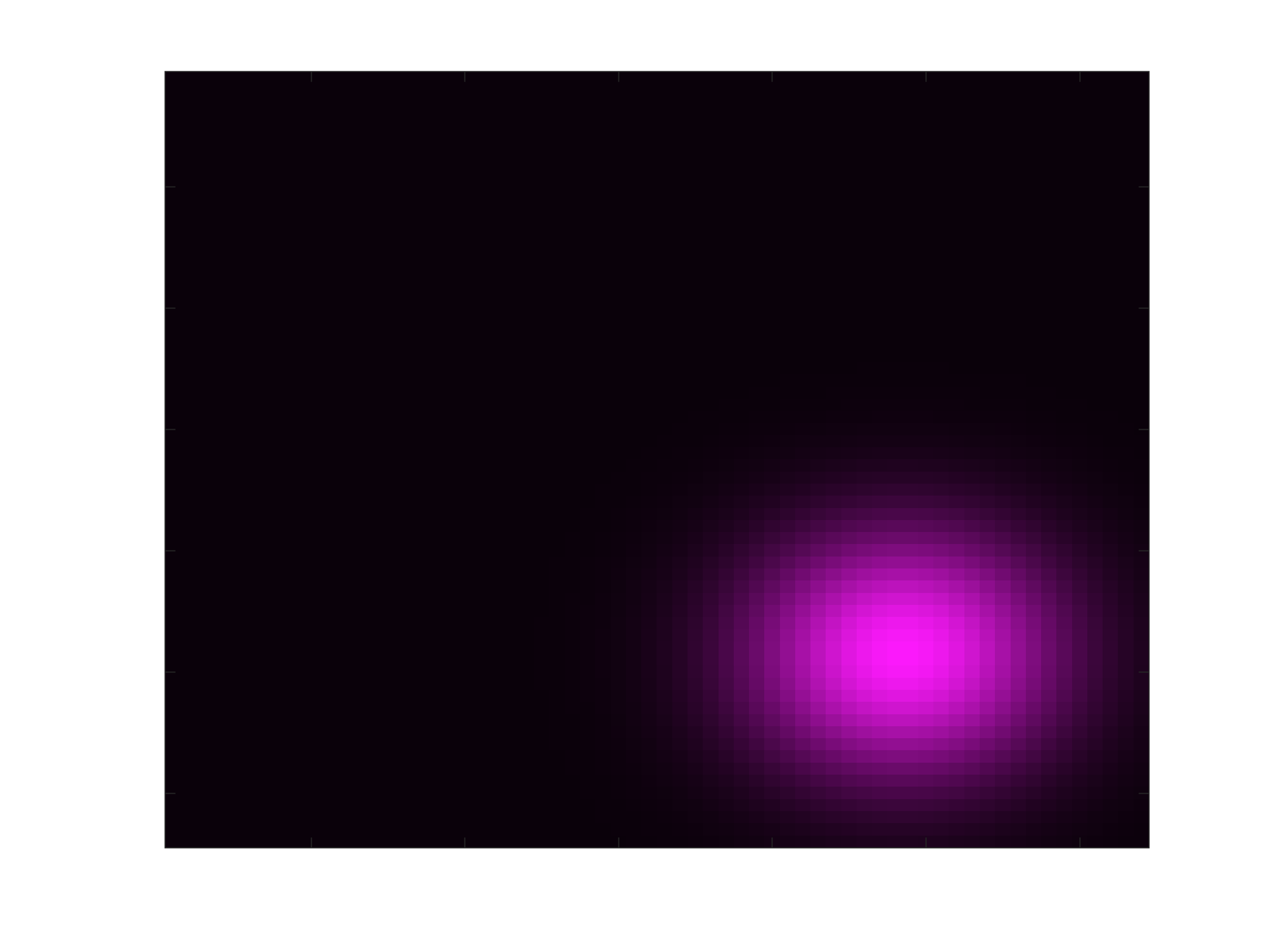}}
\subfloat[$\rho_1$]{\includegraphics[width=0.40\textwidth]{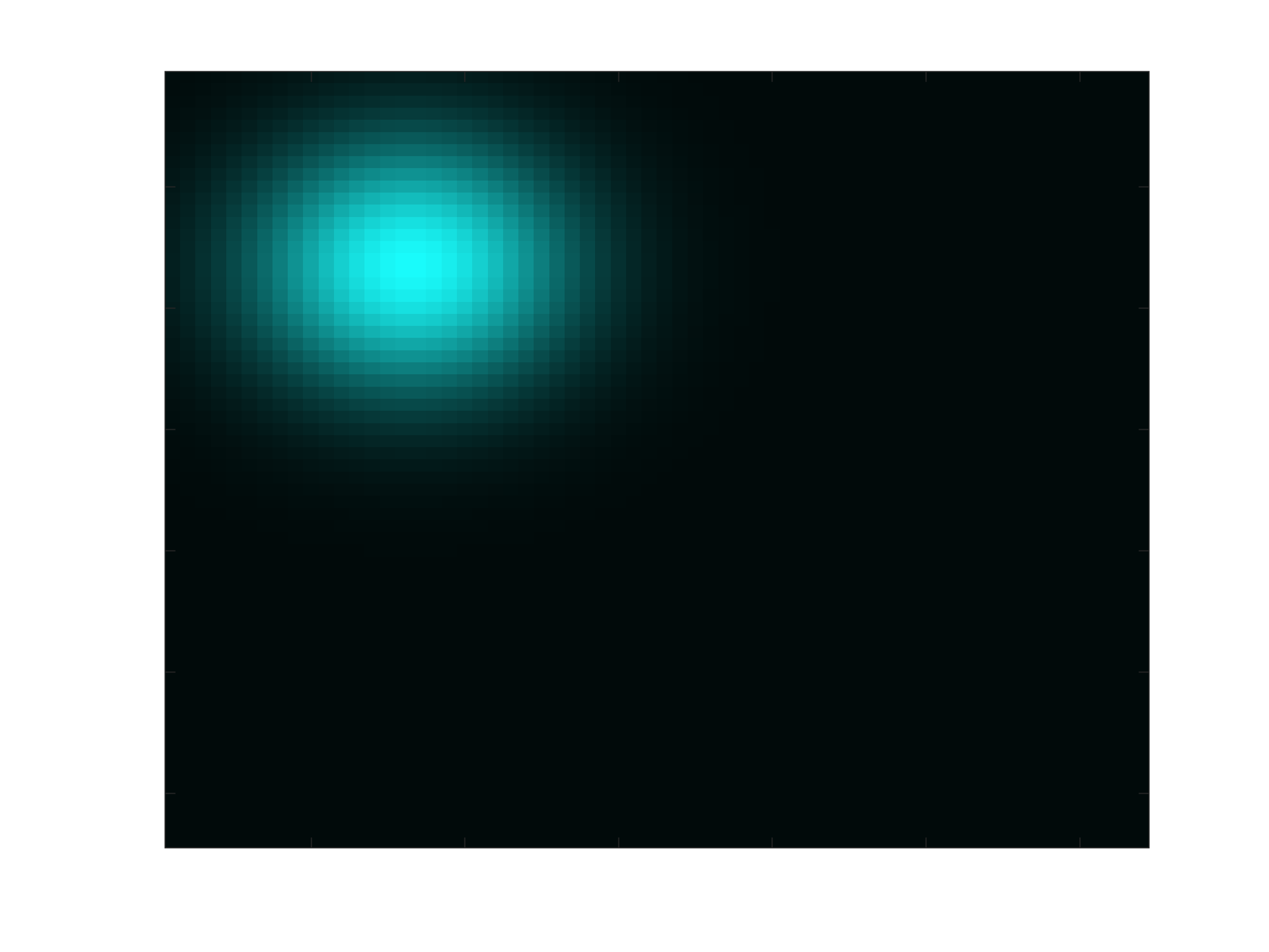}}
 \caption{Marginal distributions}
 \label{fig:ex1marginals}
\end{figure}

\begin{figure}[h]
\centering
\subfloat[$t=0.1$]{\includegraphics[width=0.2\textwidth]{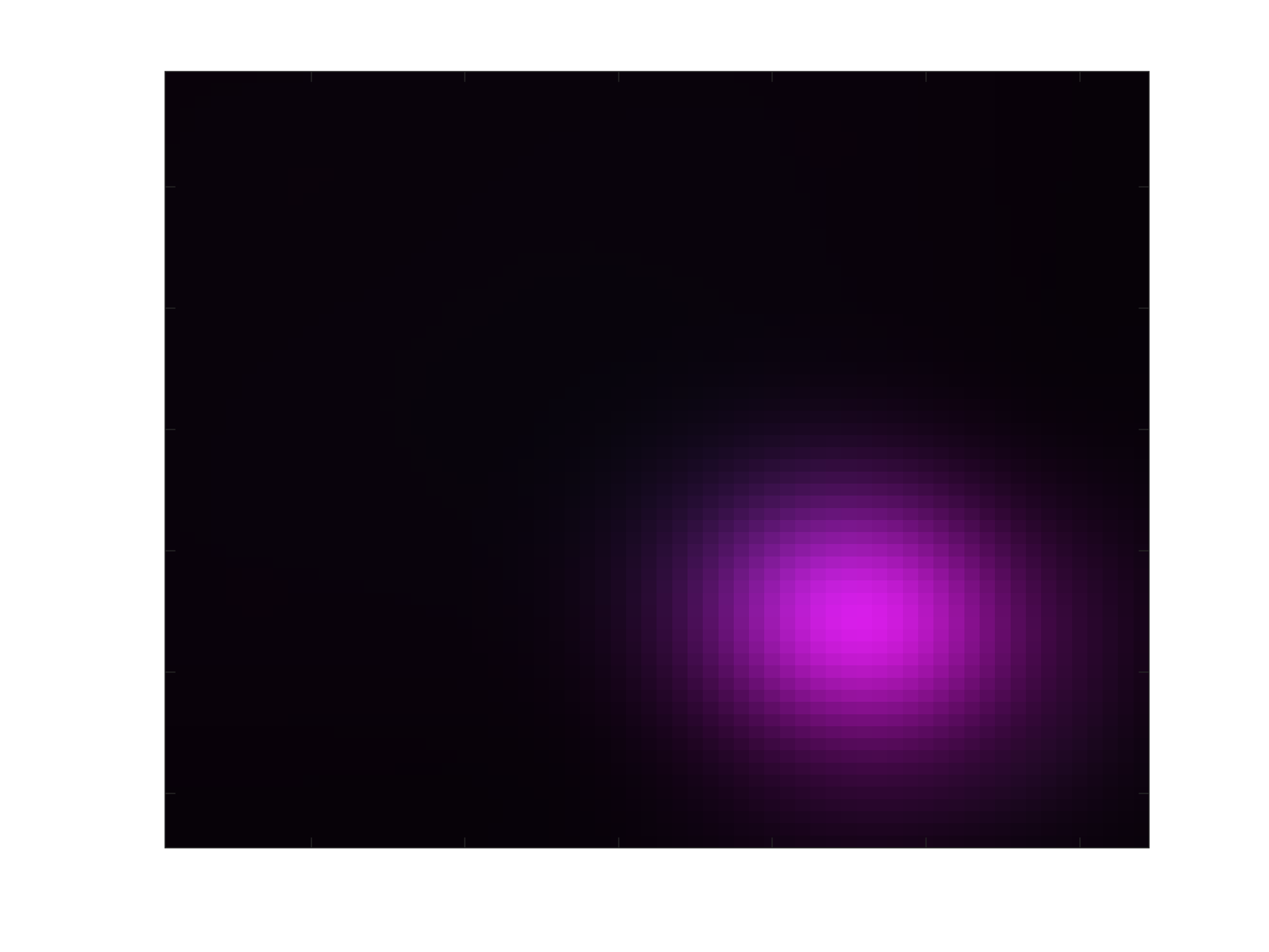}}
\subfloat[$t=0.2$]{\includegraphics[width=0.2\textwidth]{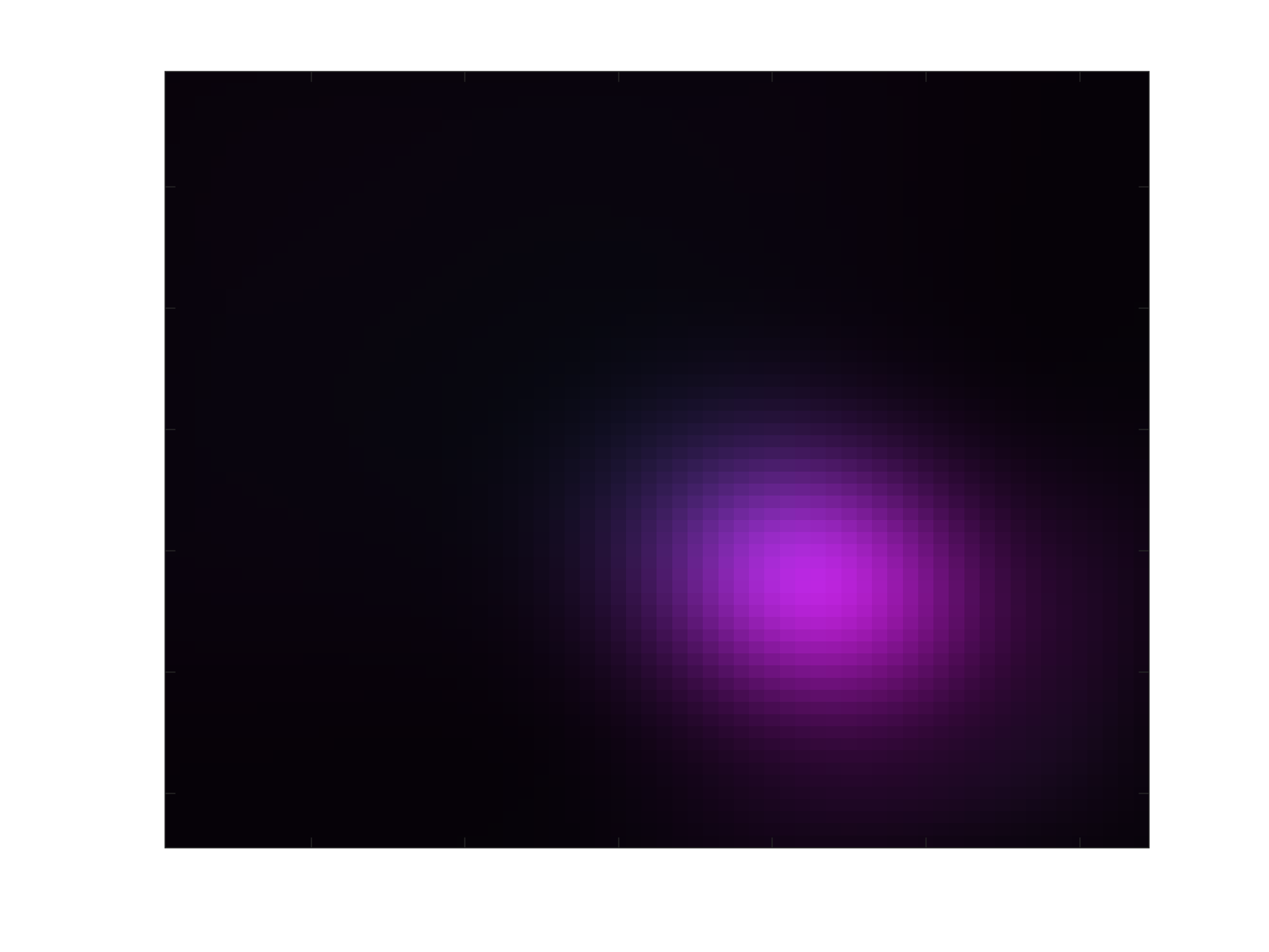}}
\subfloat[$t=0.3$]{\includegraphics[width=0.2\textwidth]{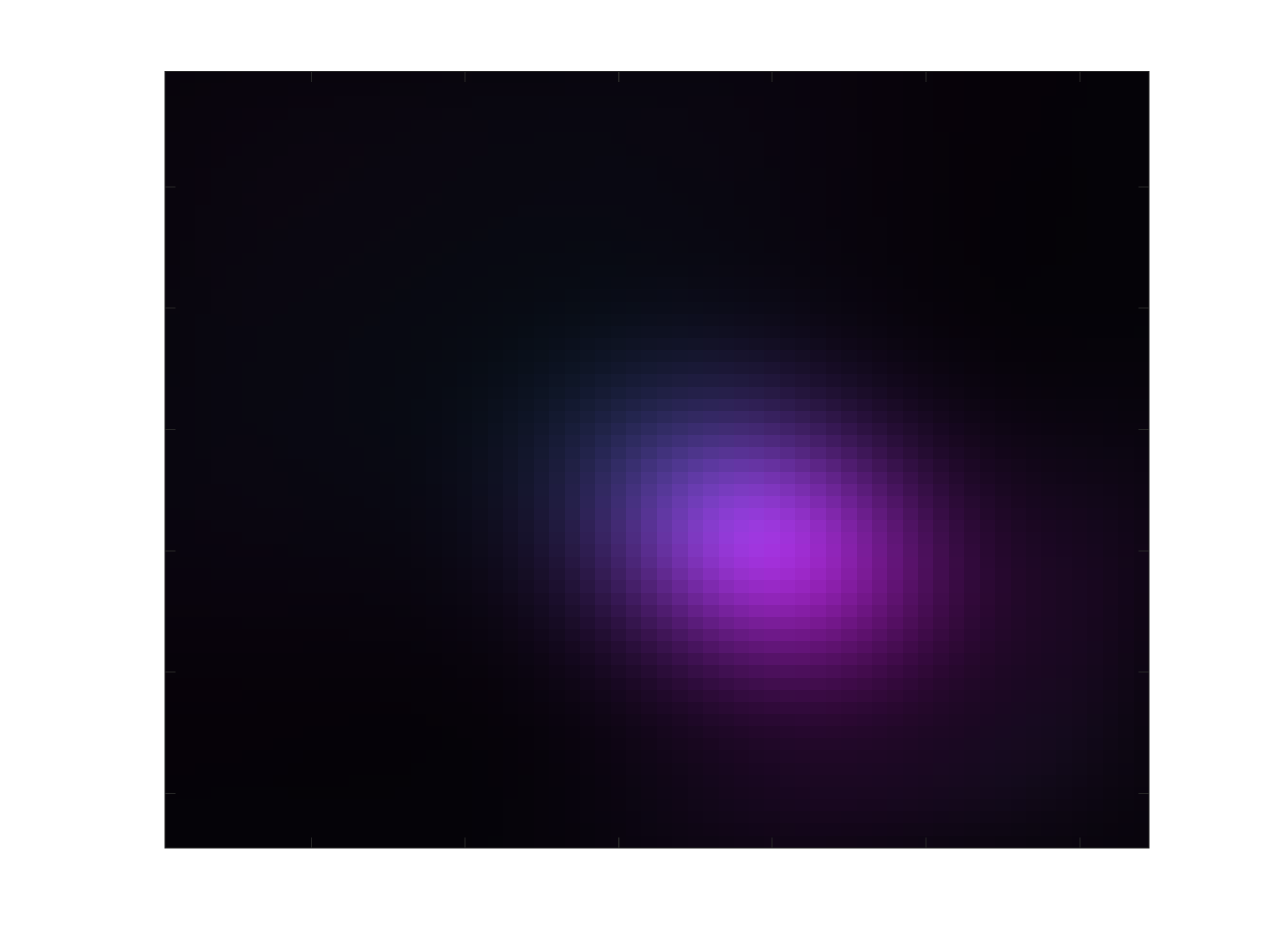}}
\subfloat[$t=0.4$]{\includegraphics[width=0.2\textwidth]{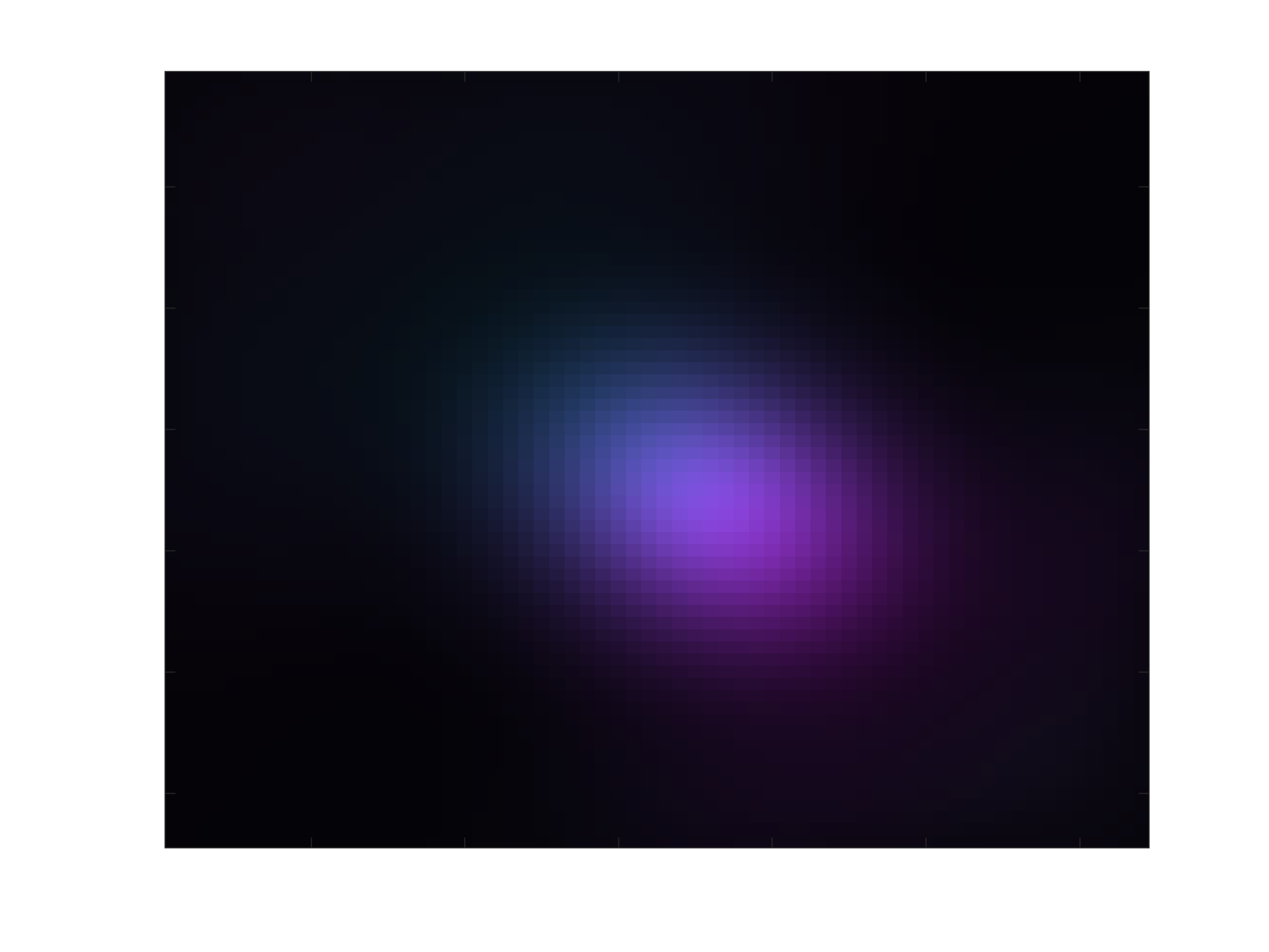}}
\\
\subfloat[$t=0.5$]{\includegraphics[width=0.2\textwidth]{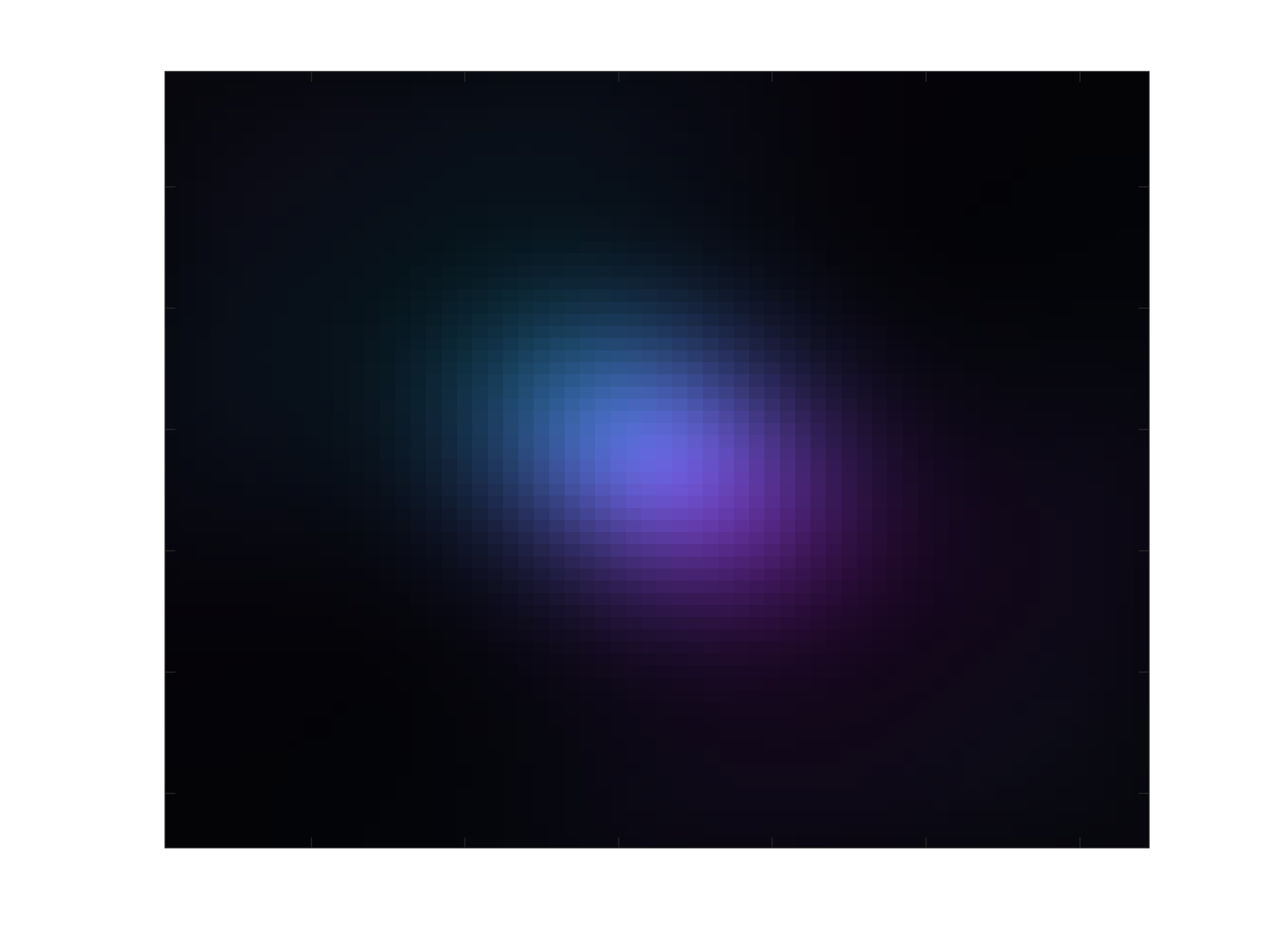}}
\subfloat[$t=0.6$]{\includegraphics[width=0.2\textwidth]{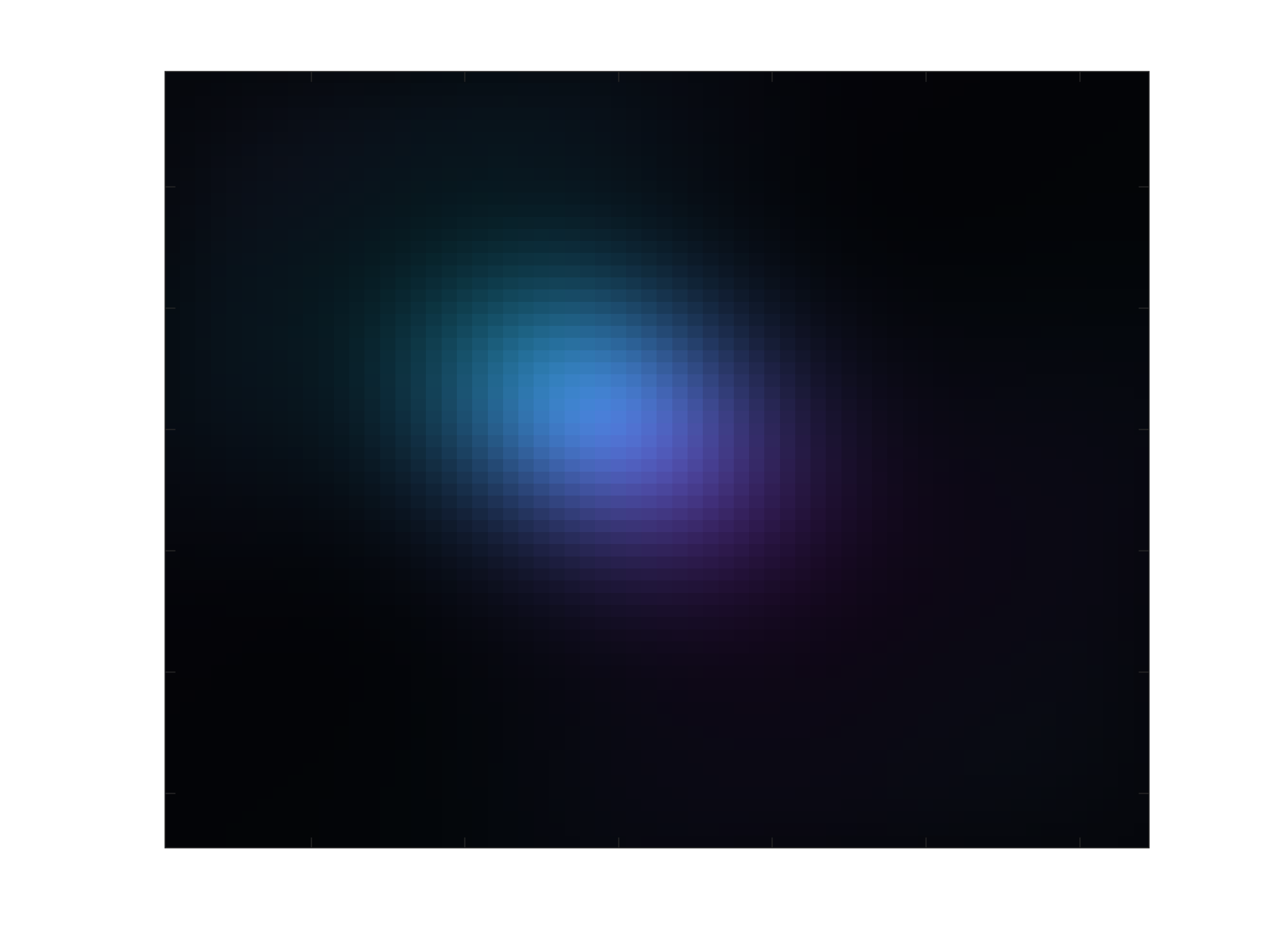}}
\subfloat[$t=0.7$]{\includegraphics[width=0.2\textwidth]{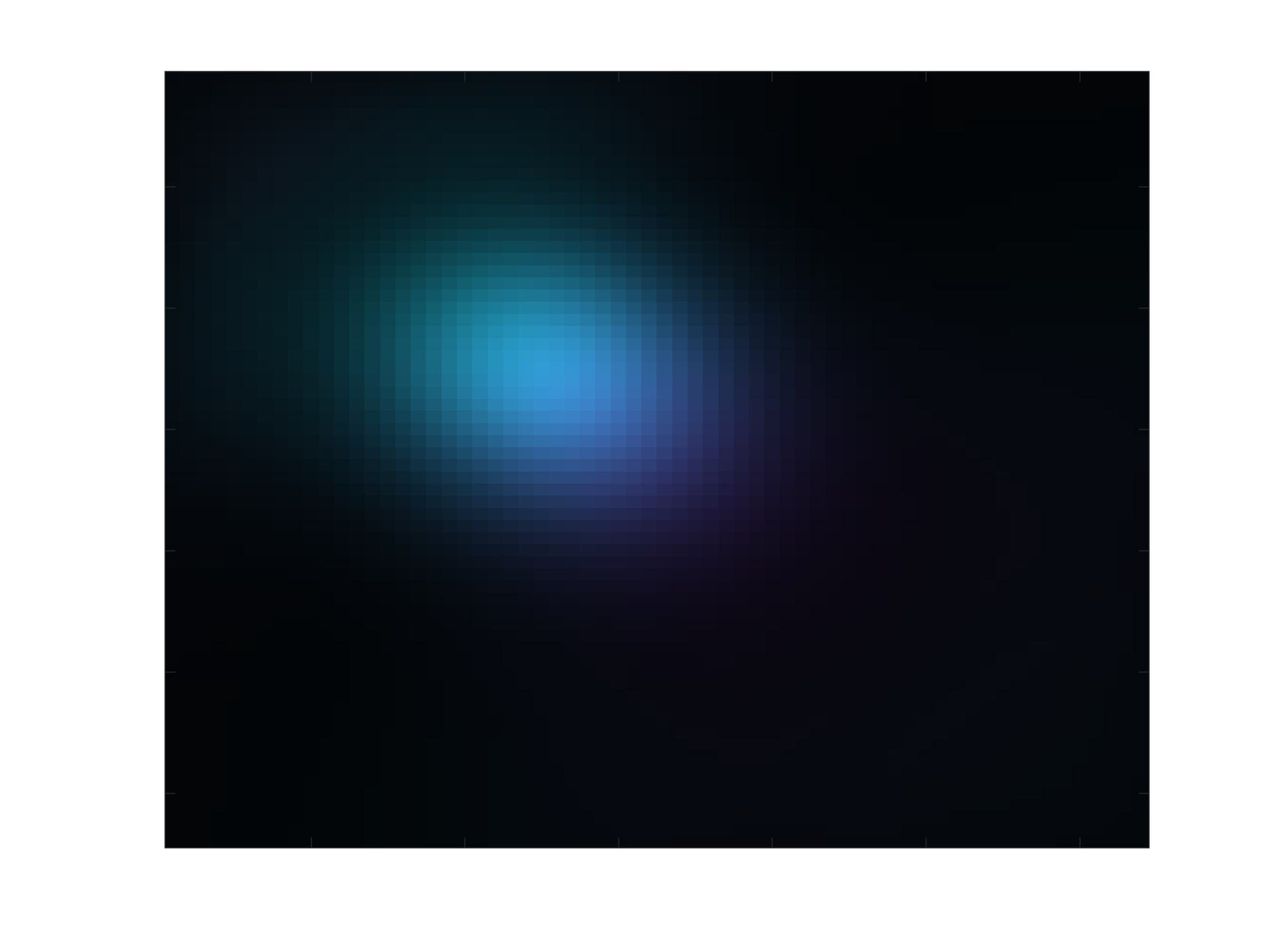}}
\subfloat[$t=0.8$]{\includegraphics[width=0.2\textwidth]{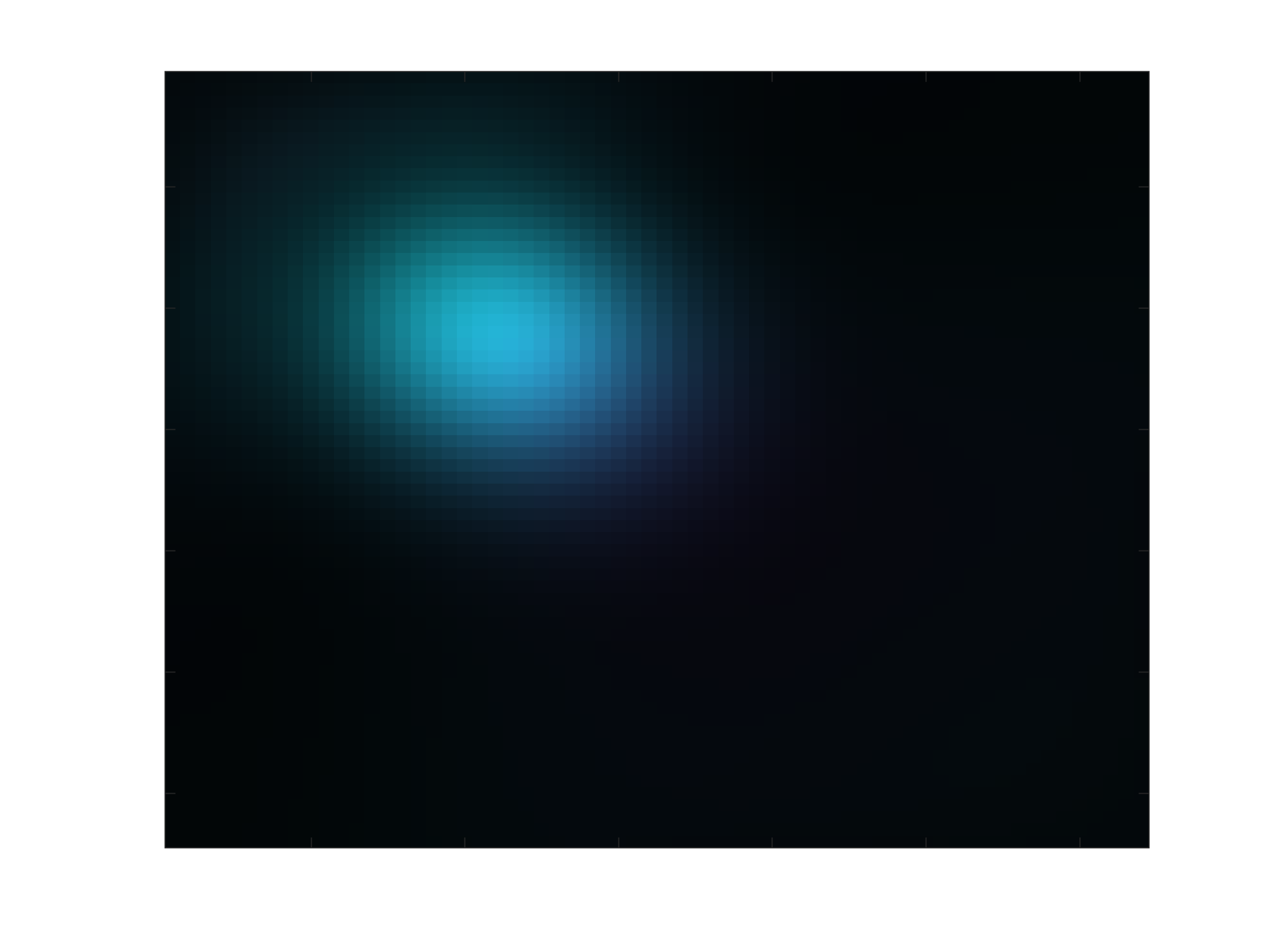}}
\subfloat[$t=0.9$]{\includegraphics[width=0.2\textwidth]{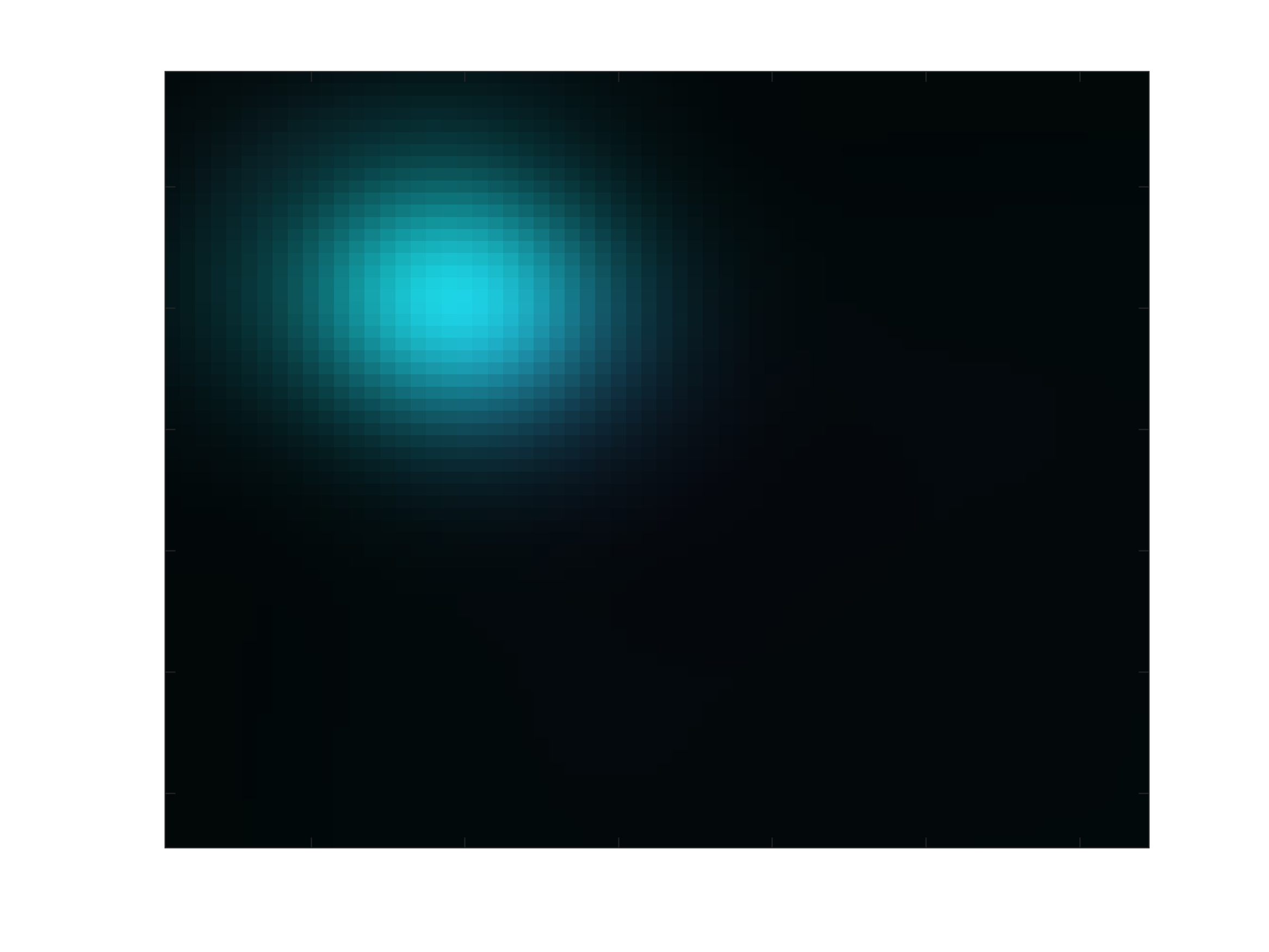}}
\caption{Interpolation with $\gamma=0.001$}
\label{fig:ex1interp}
\end{figure}

\begin{figure}[h]
\centering
\subfloat[$\rho_0$]{\includegraphics[width=0.40\textwidth]{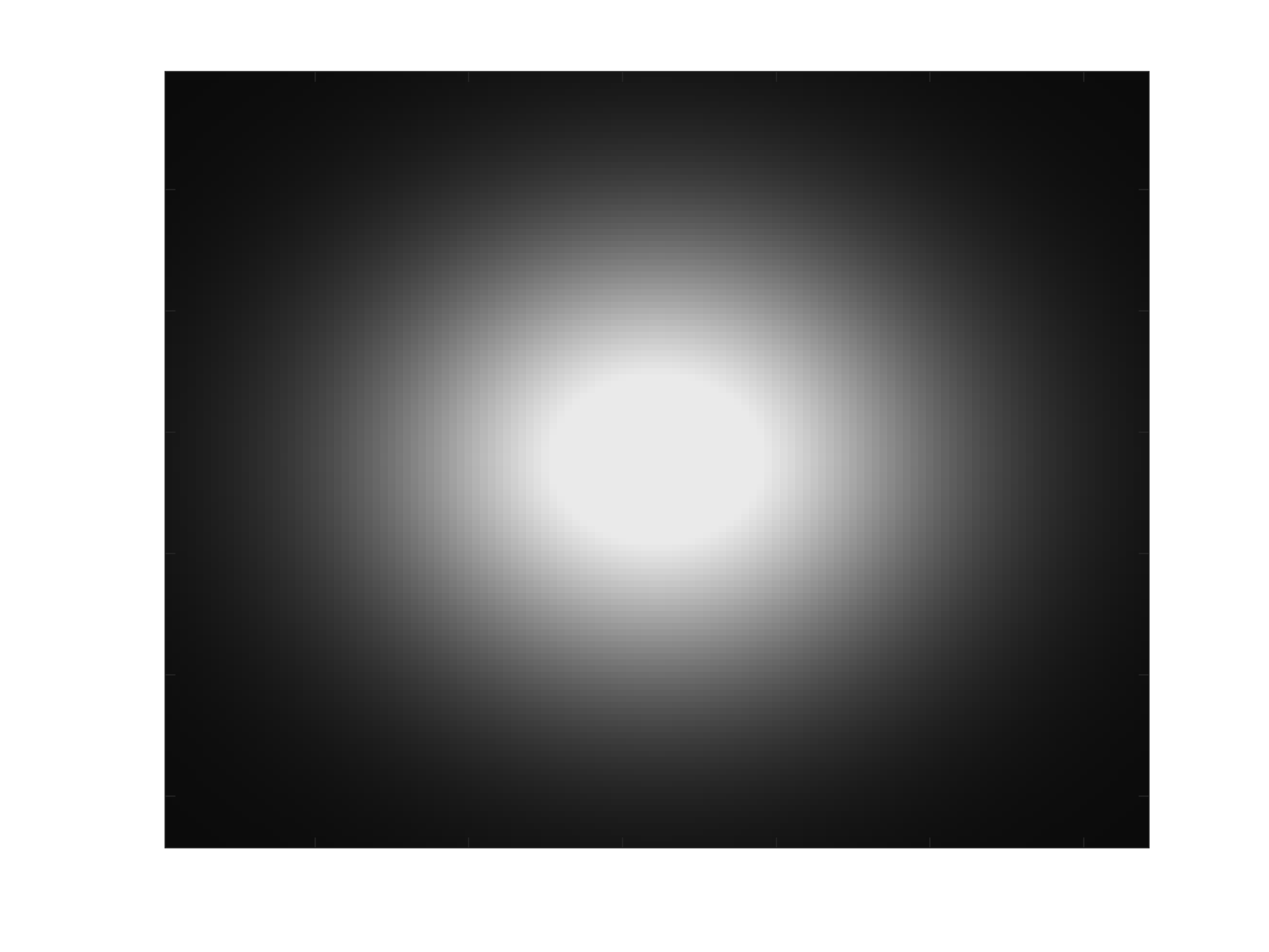}}
\subfloat[$\rho_1$]{\includegraphics[width=0.40\textwidth]{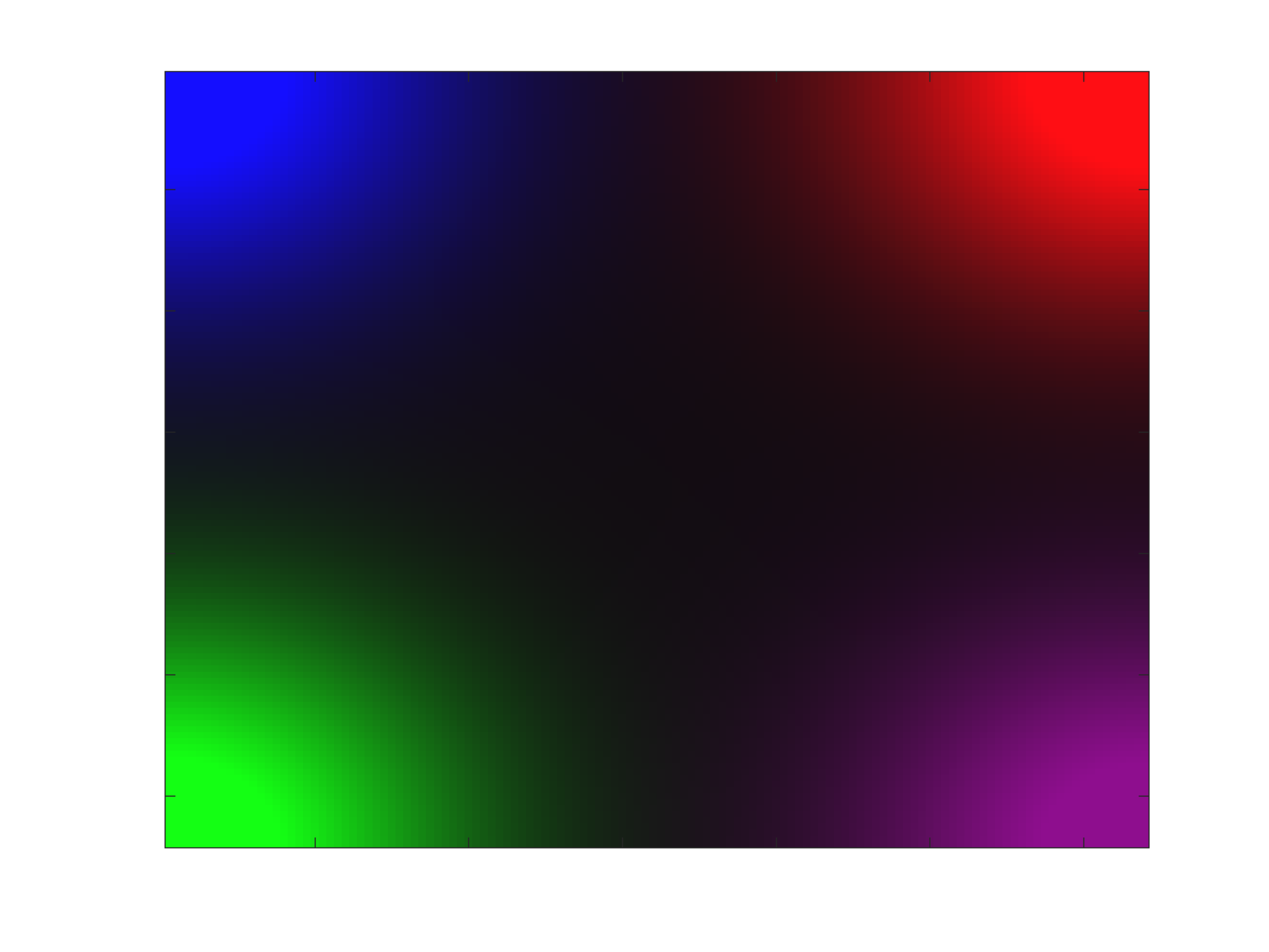}}
 \caption{Marginal distributions}
 \label{fig:ex2marginals}
\end{figure}

\begin{figure}[h]
\centering
\subfloat[$t=0.1$]{\includegraphics[width=0.2\textwidth]{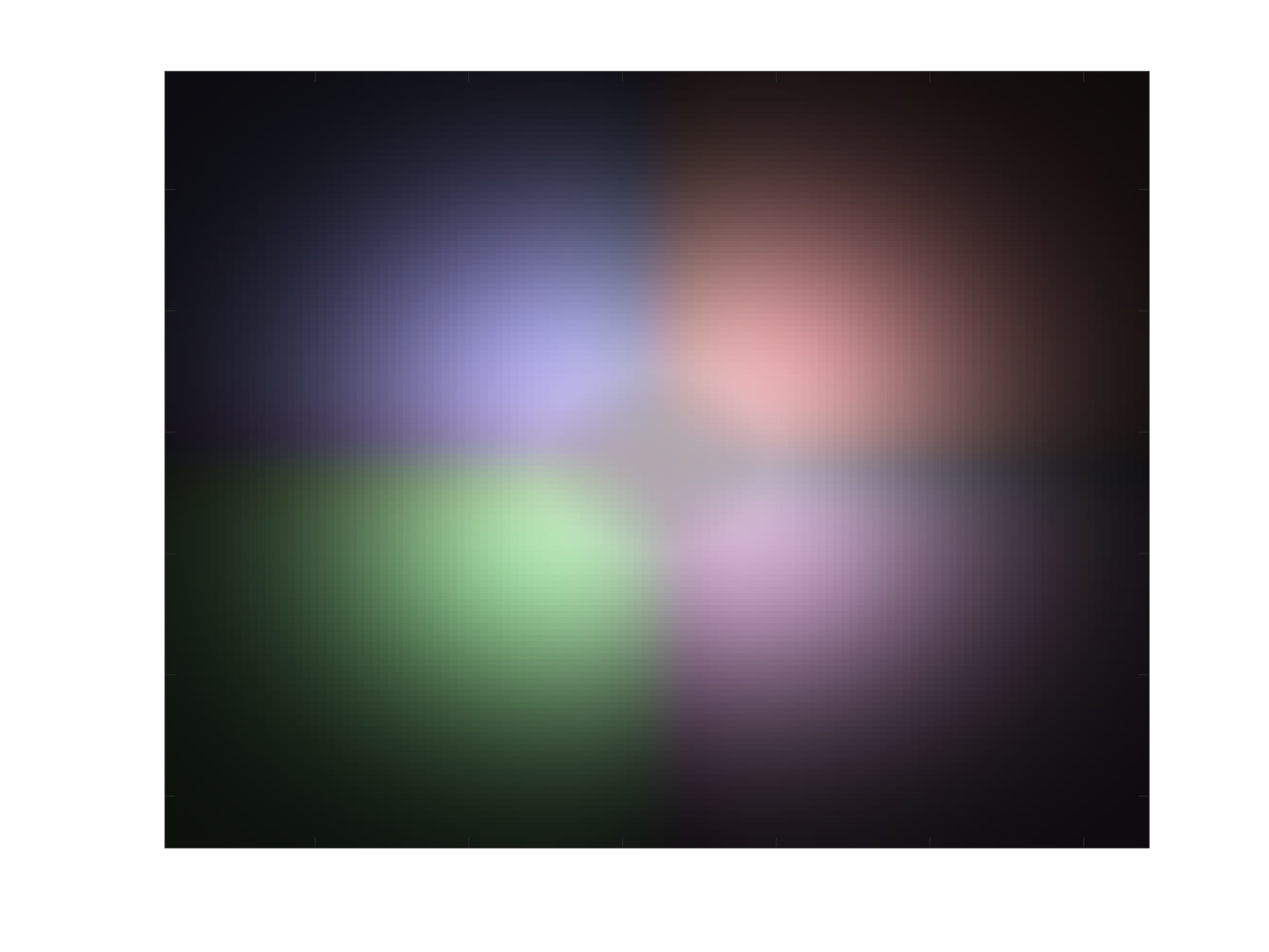}}
\subfloat[$t=0.2$]{\includegraphics[width=0.2\textwidth]{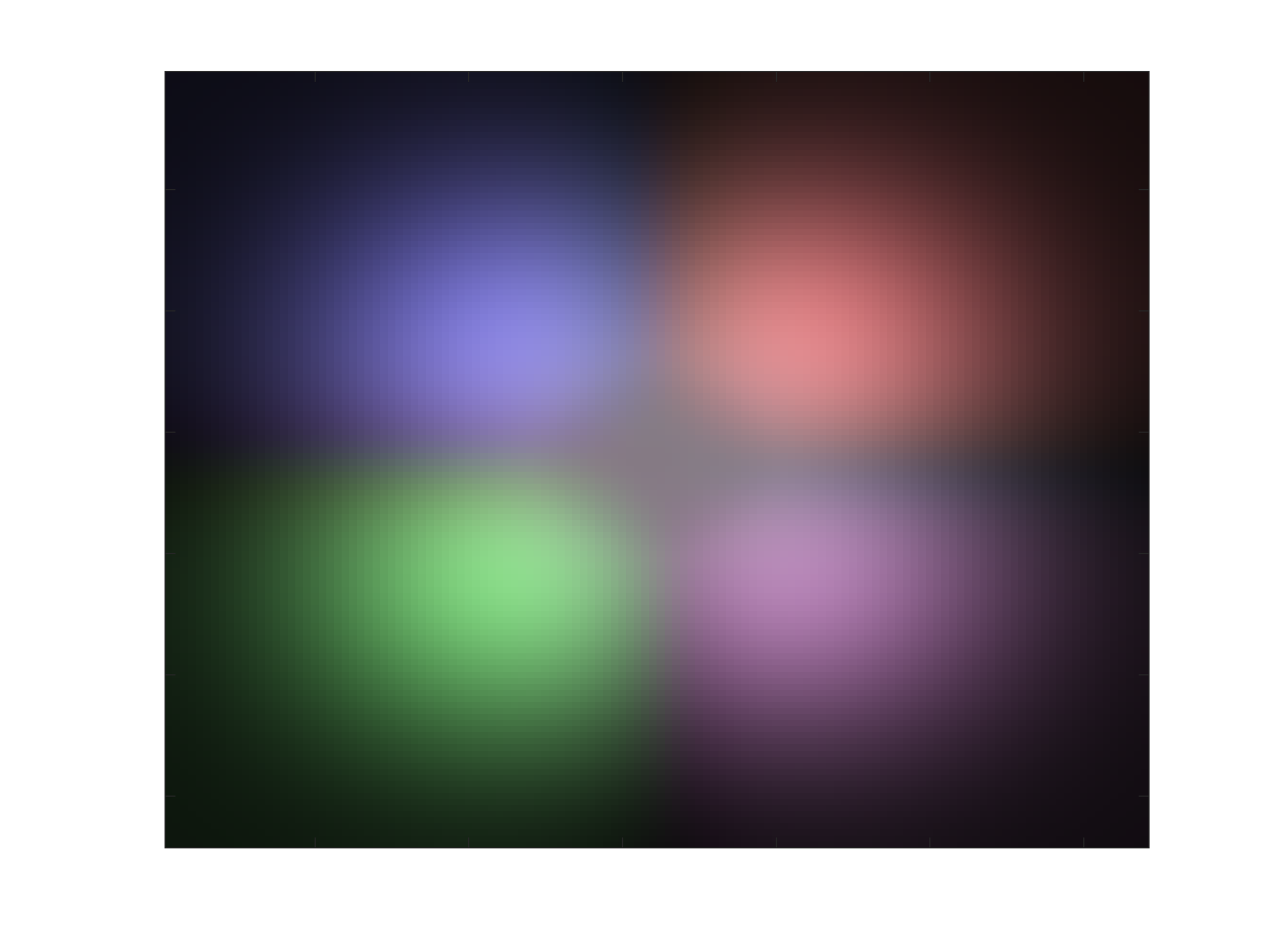}}
\subfloat[$t=0.3$]{\includegraphics[width=0.2\textwidth]{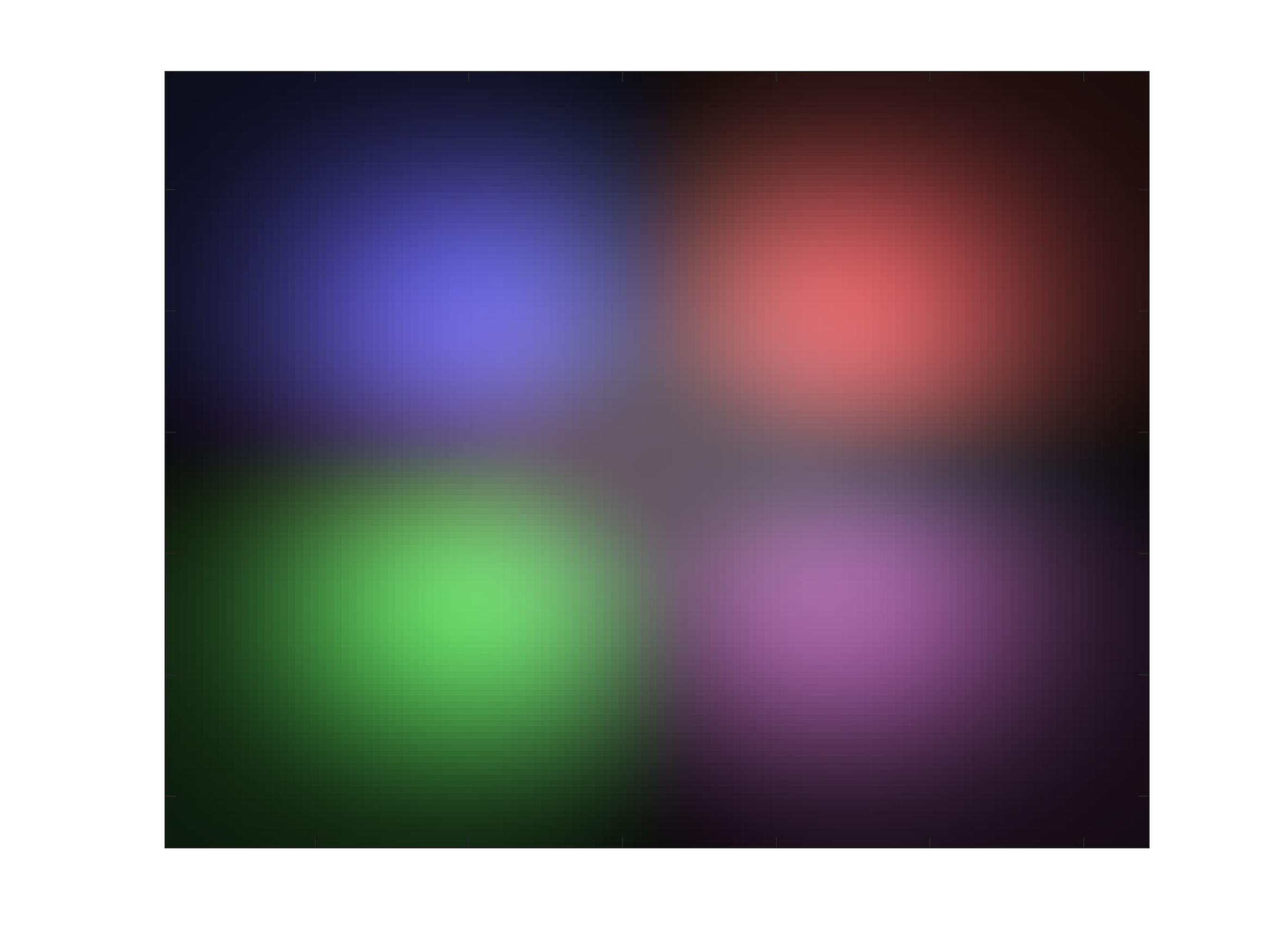}}
\subfloat[$t=0.4$]{\includegraphics[width=0.2\textwidth]{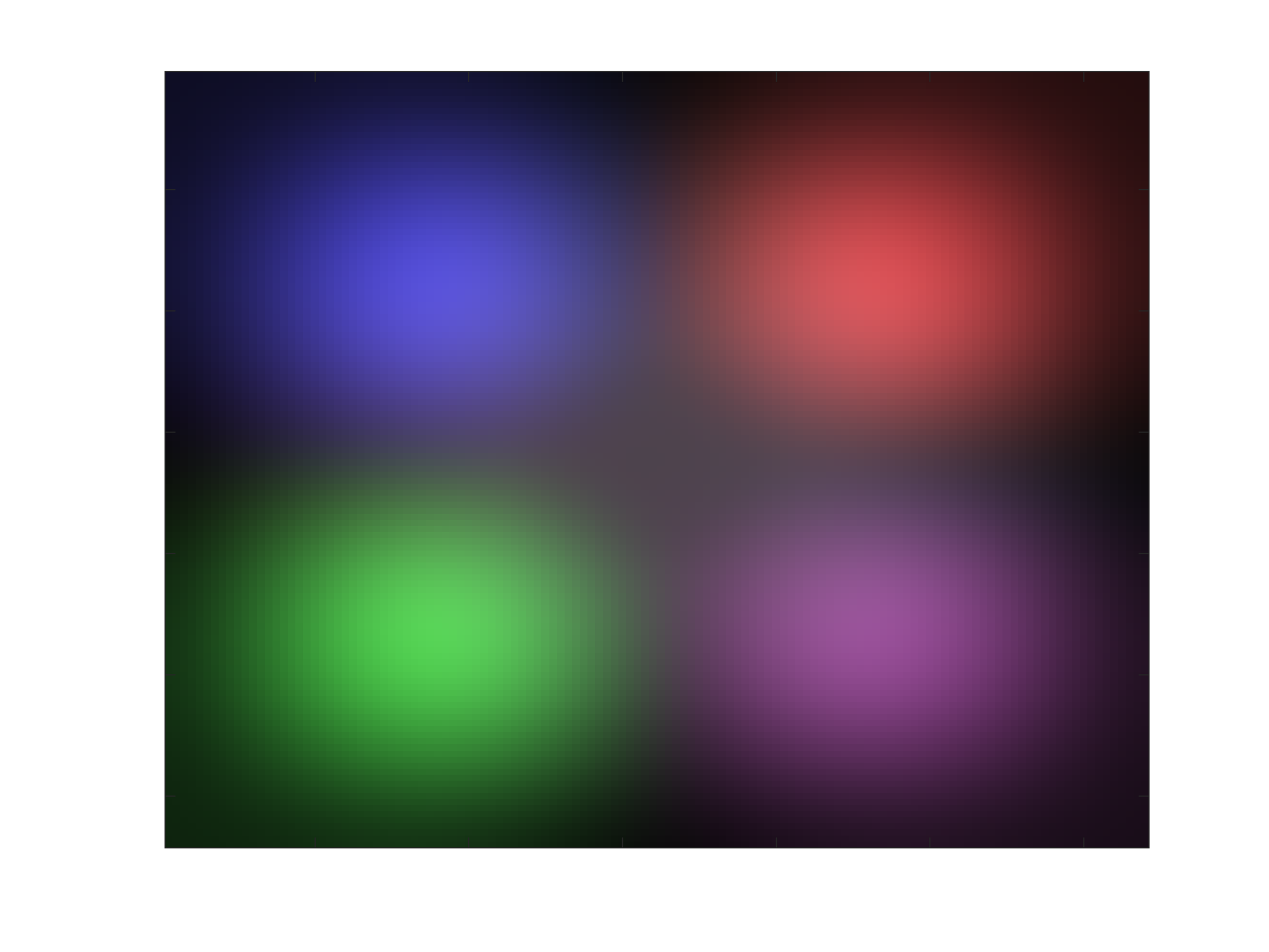}}
\\
\subfloat[$t=0.5$]{\includegraphics[width=0.2\textwidth]{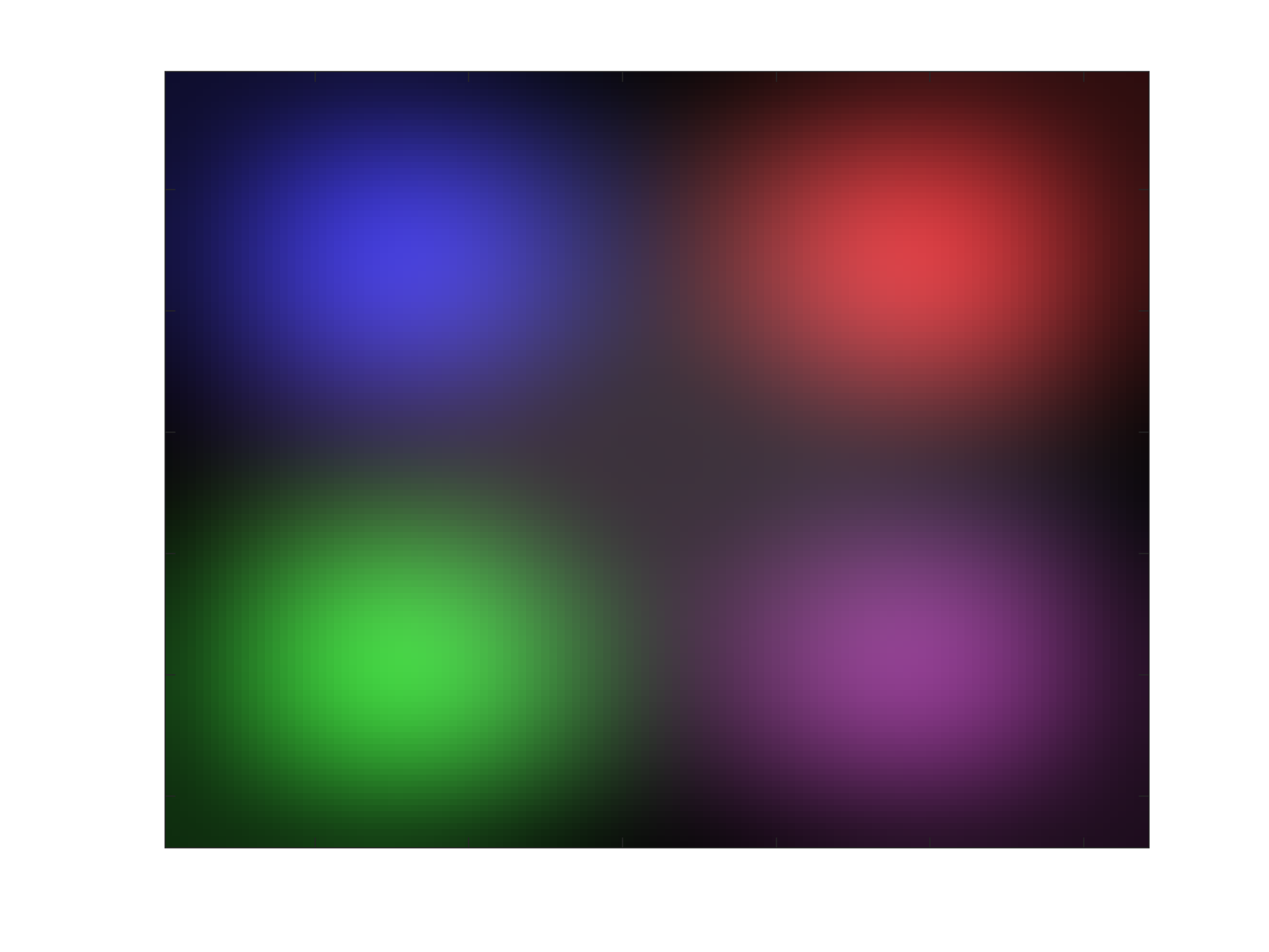}}
\subfloat[$t=0.6$]{\includegraphics[width=0.2\textwidth]{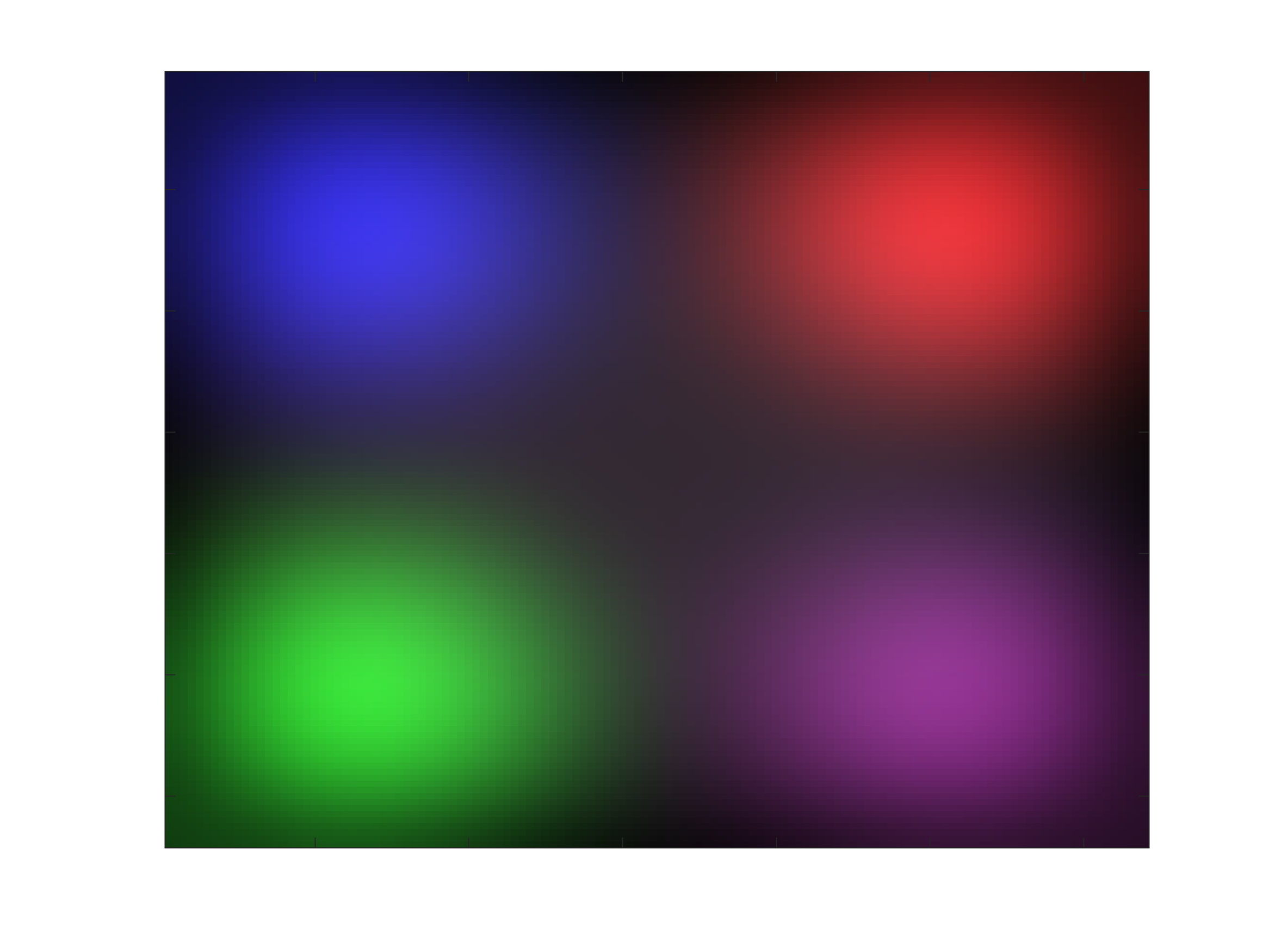}}
\subfloat[$t=0.7$]{\includegraphics[width=0.2\textwidth]{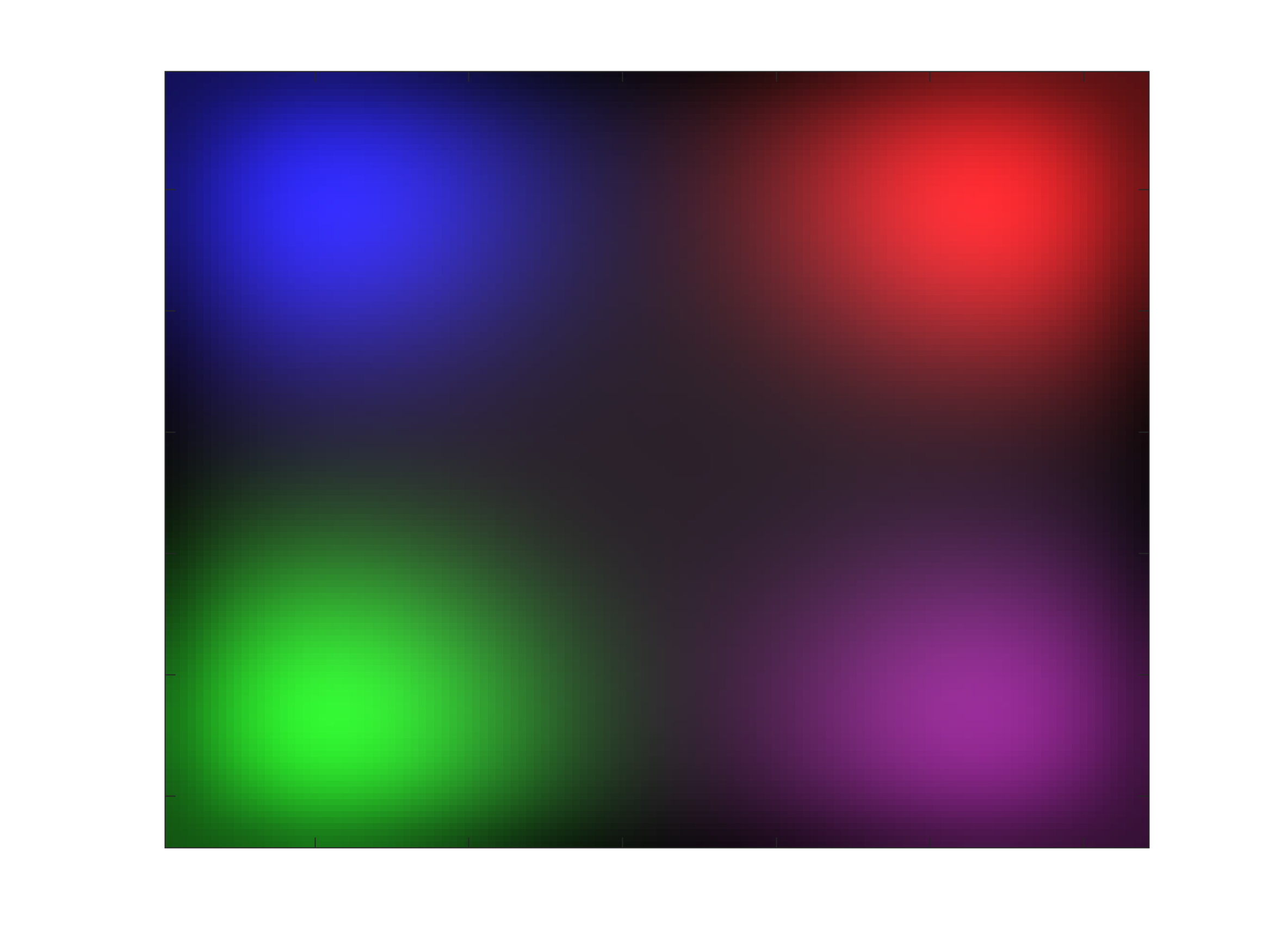}}
\subfloat[$t=0.8$]{\includegraphics[width=0.2\textwidth]{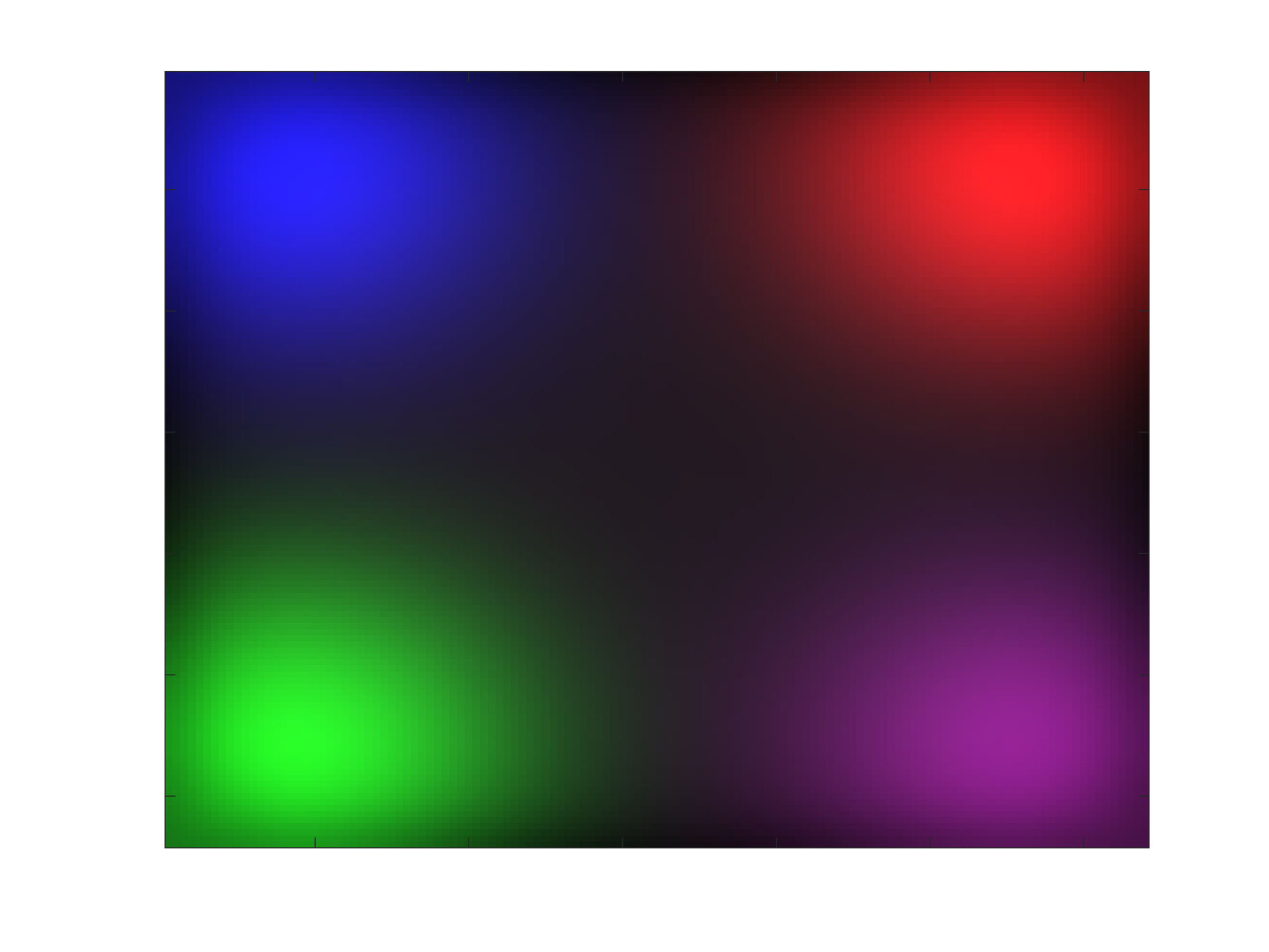}}
\subfloat[$t=0.9$]{\includegraphics[width=0.2\textwidth]{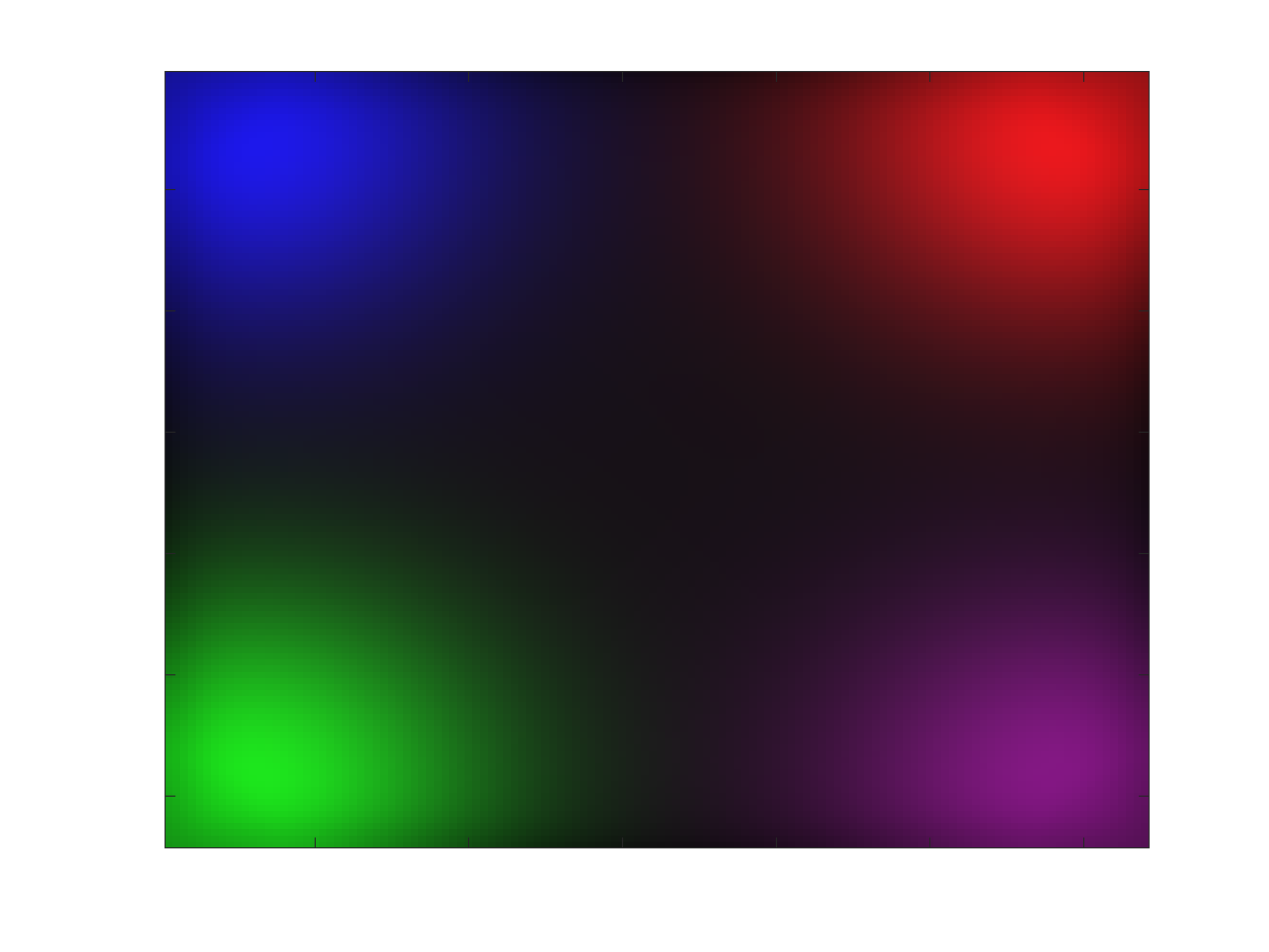}}
\caption{Interpolation with $\gamma=0.01$}
\label{fig:ex2interp}
\end{figure}

\begin{figure}[h]
\centering
\subfloat[$\rho_0$]{\includegraphics[width=0.40\textwidth]{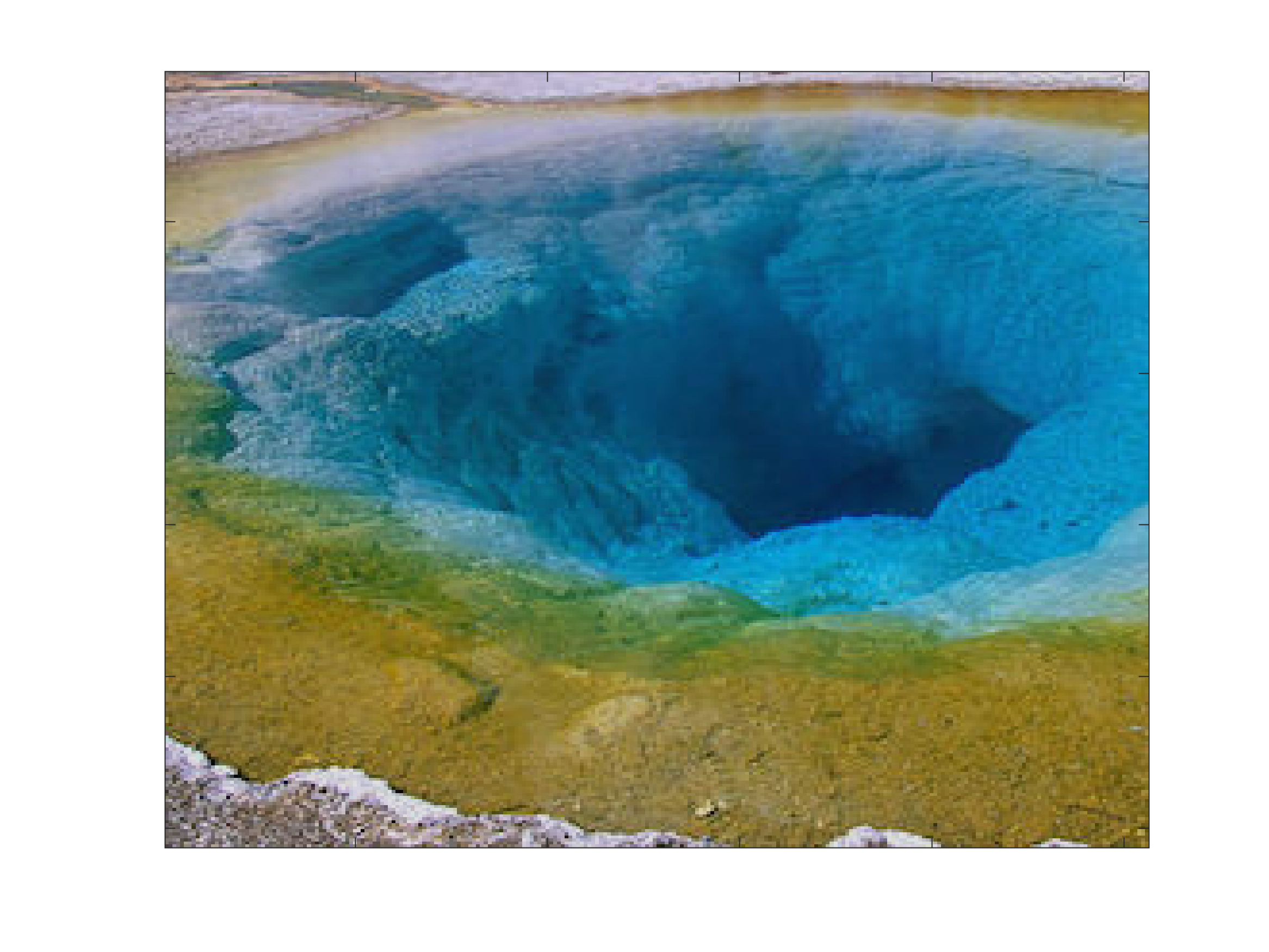}}
\subfloat[$\rho_1$]{\includegraphics[width=0.40\textwidth]{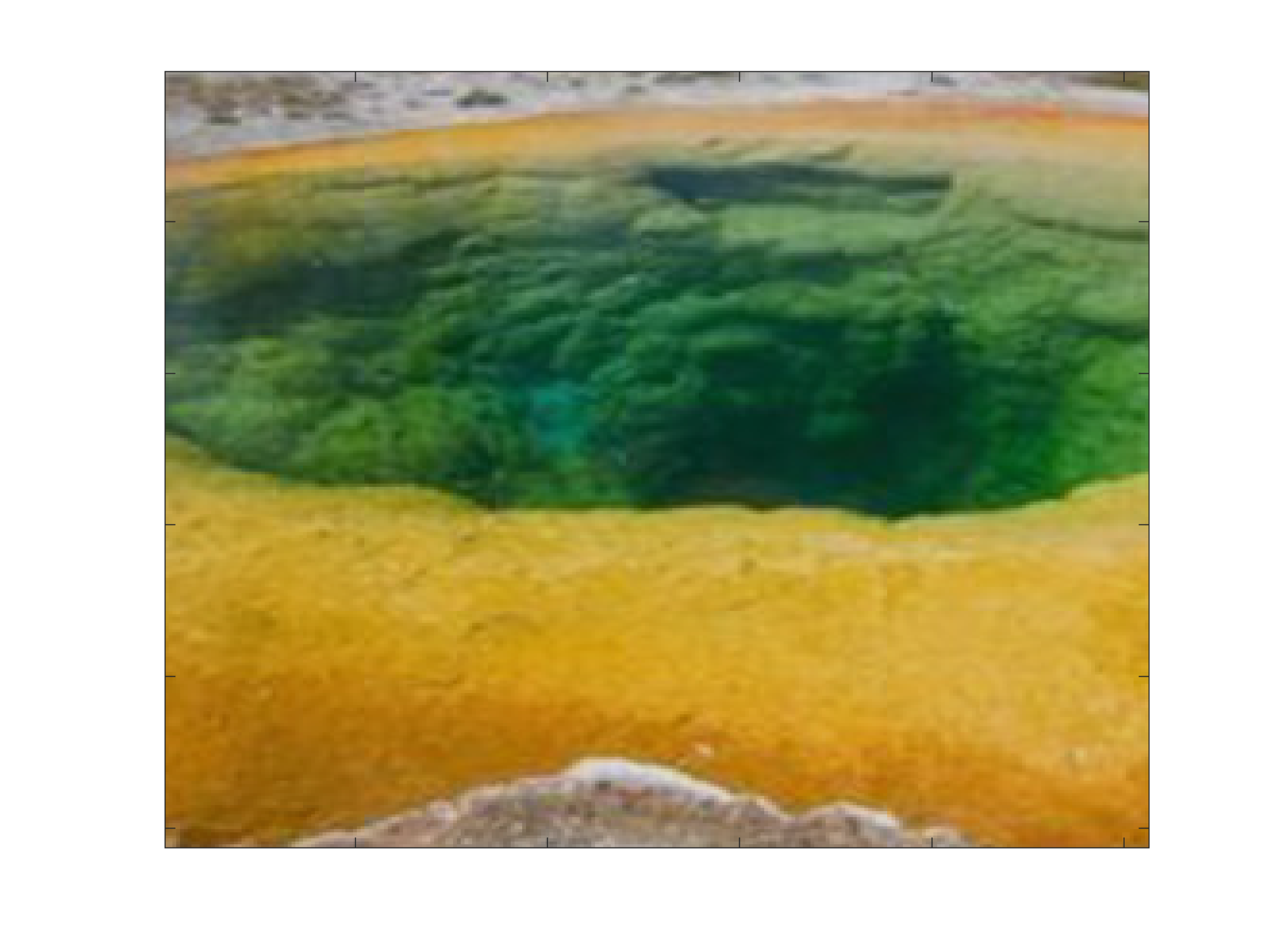}}
 \caption{Marginal distributions}
 \label{fig:ex3marginals}
\end{figure}

\begin{figure}[h]
\centering
\subfloat[$t=0.1$]{\includegraphics[width=0.35\textwidth]{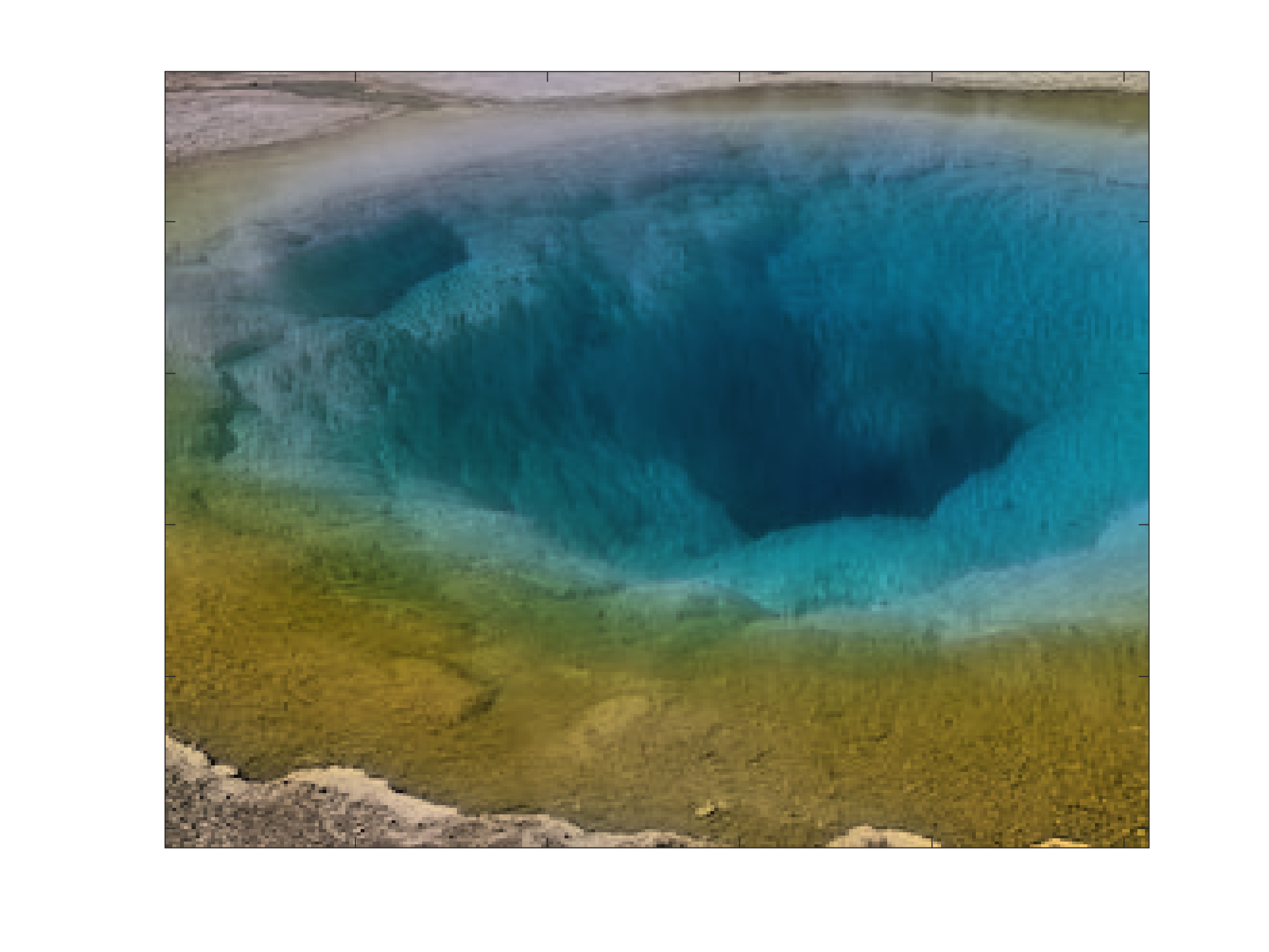}}
\subfloat[$t=0.2$]{\includegraphics[width=0.35\textwidth]{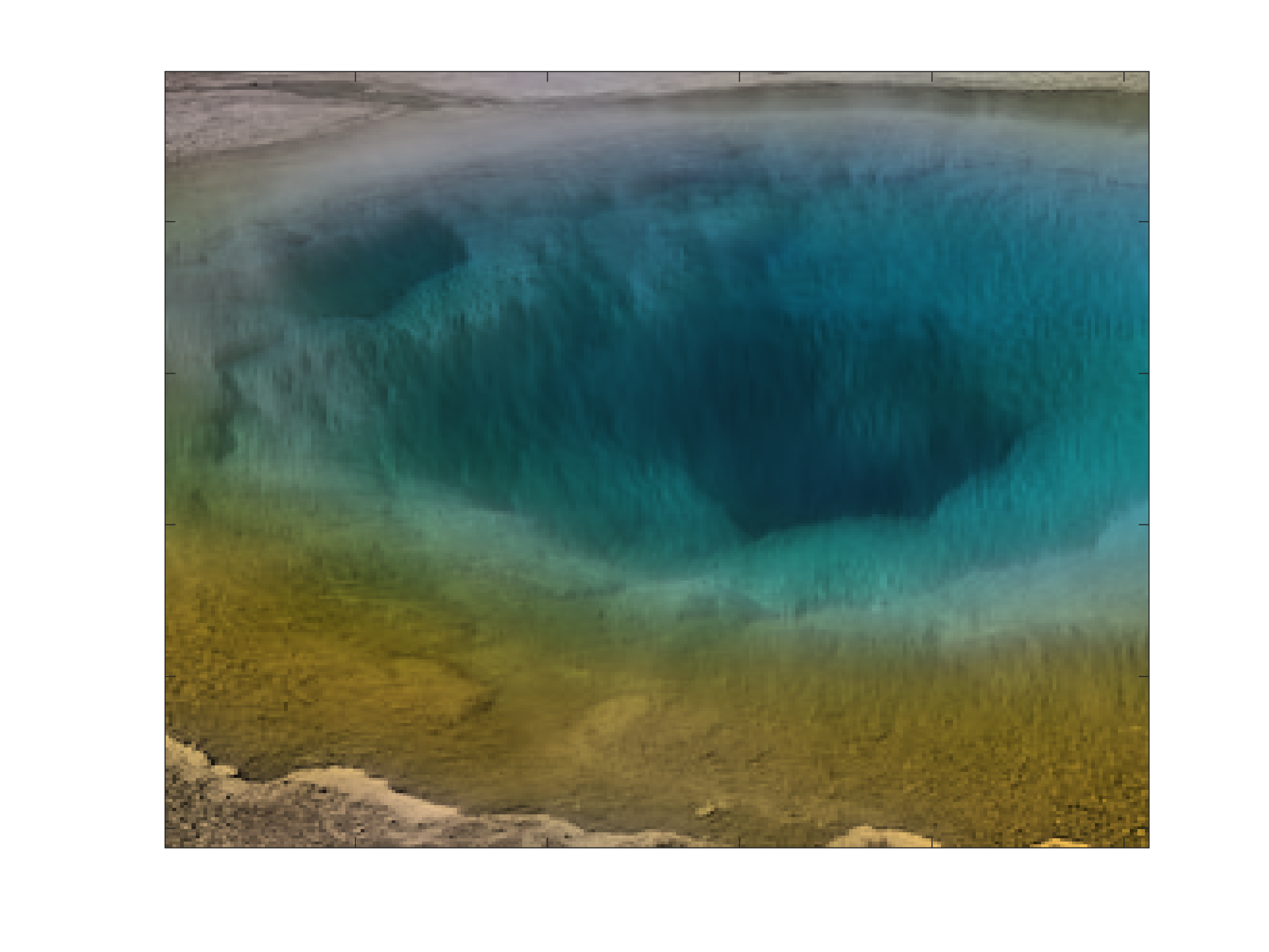}}
\subfloat[$t=0.3$]{\includegraphics[width=0.35\textwidth]{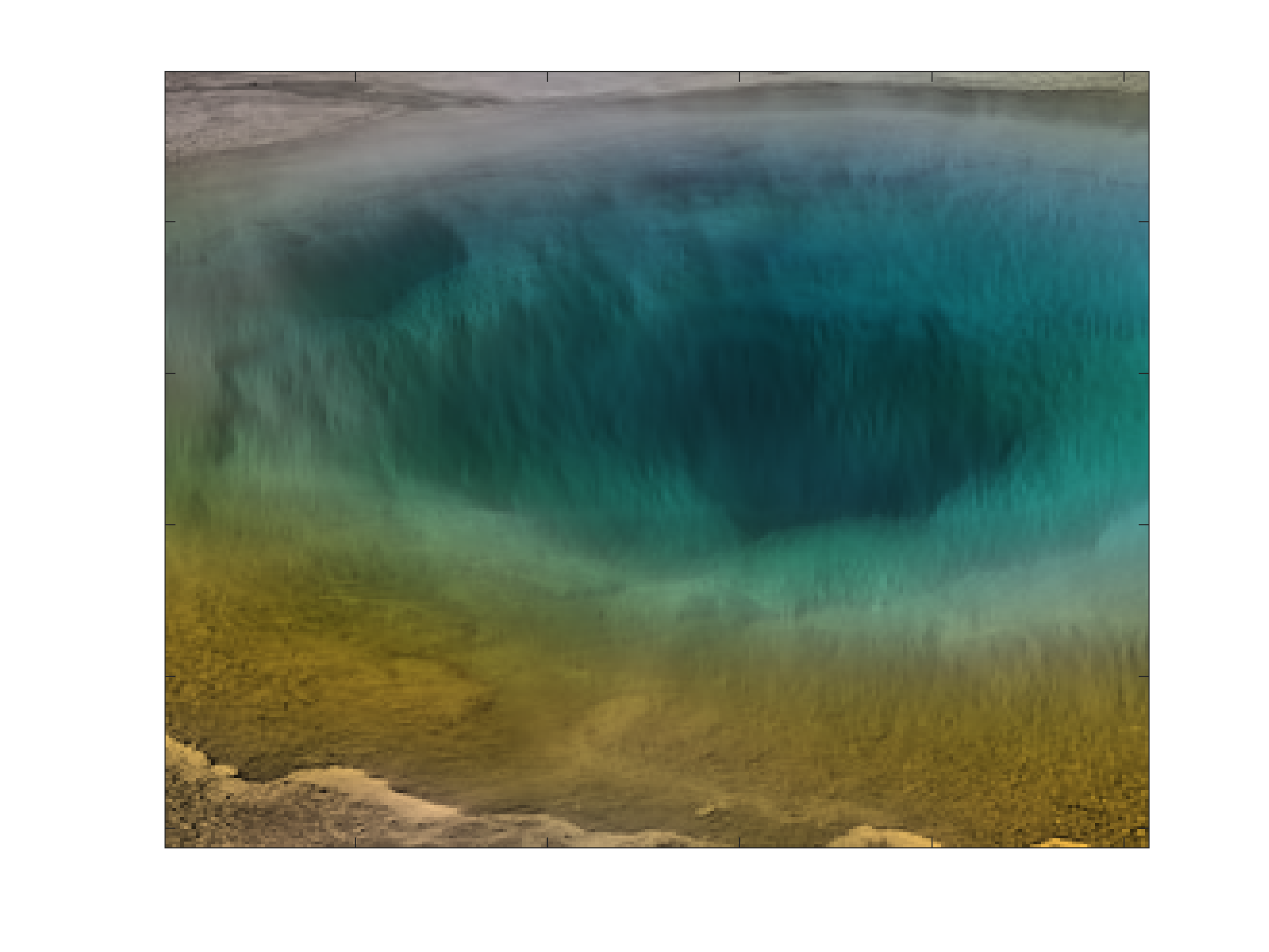}}\\
\subfloat[$t=0.4$]{\includegraphics[width=0.35\textwidth]{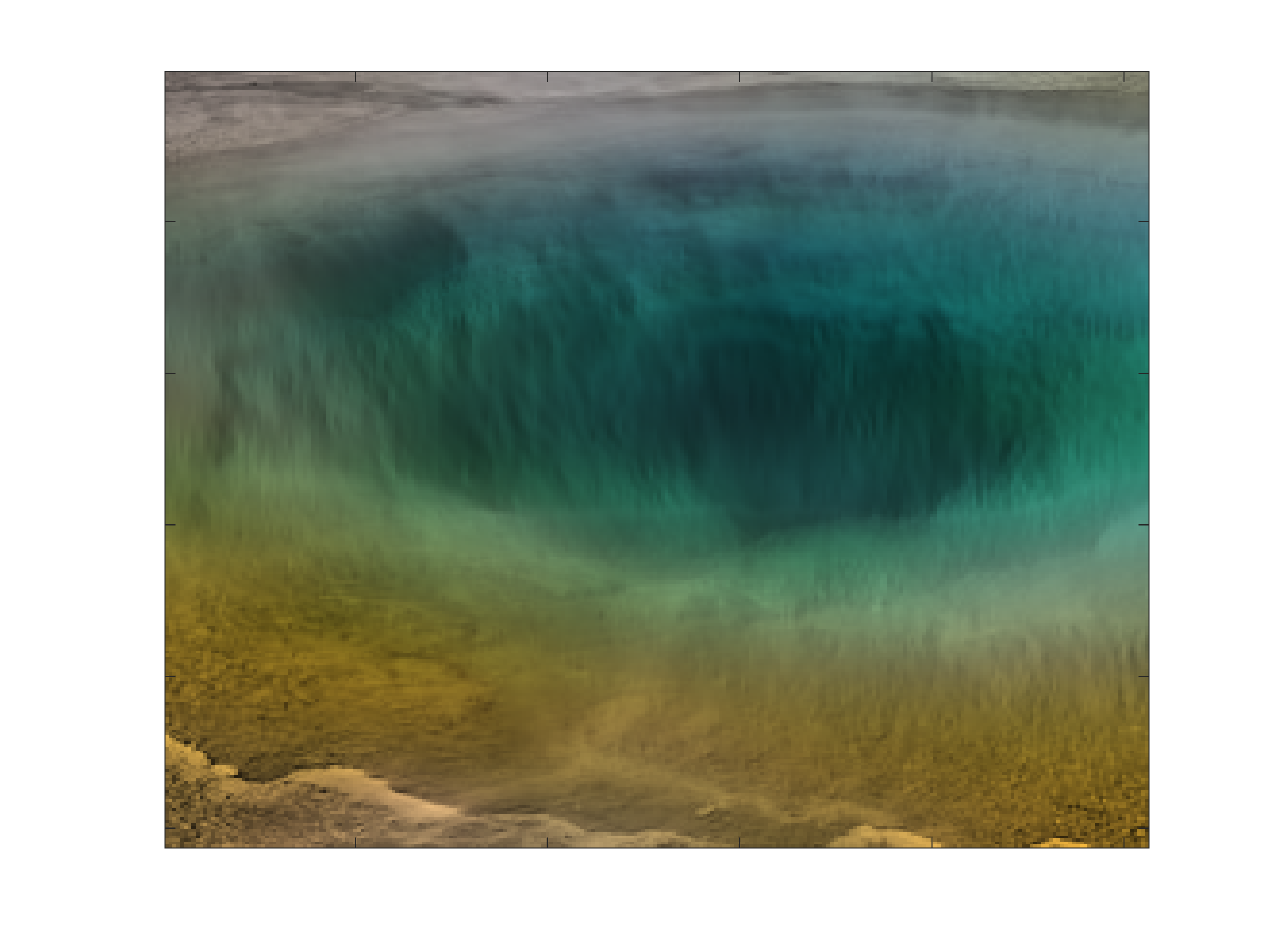}}
\subfloat[$t=0.5$]{\includegraphics[width=0.35\textwidth]{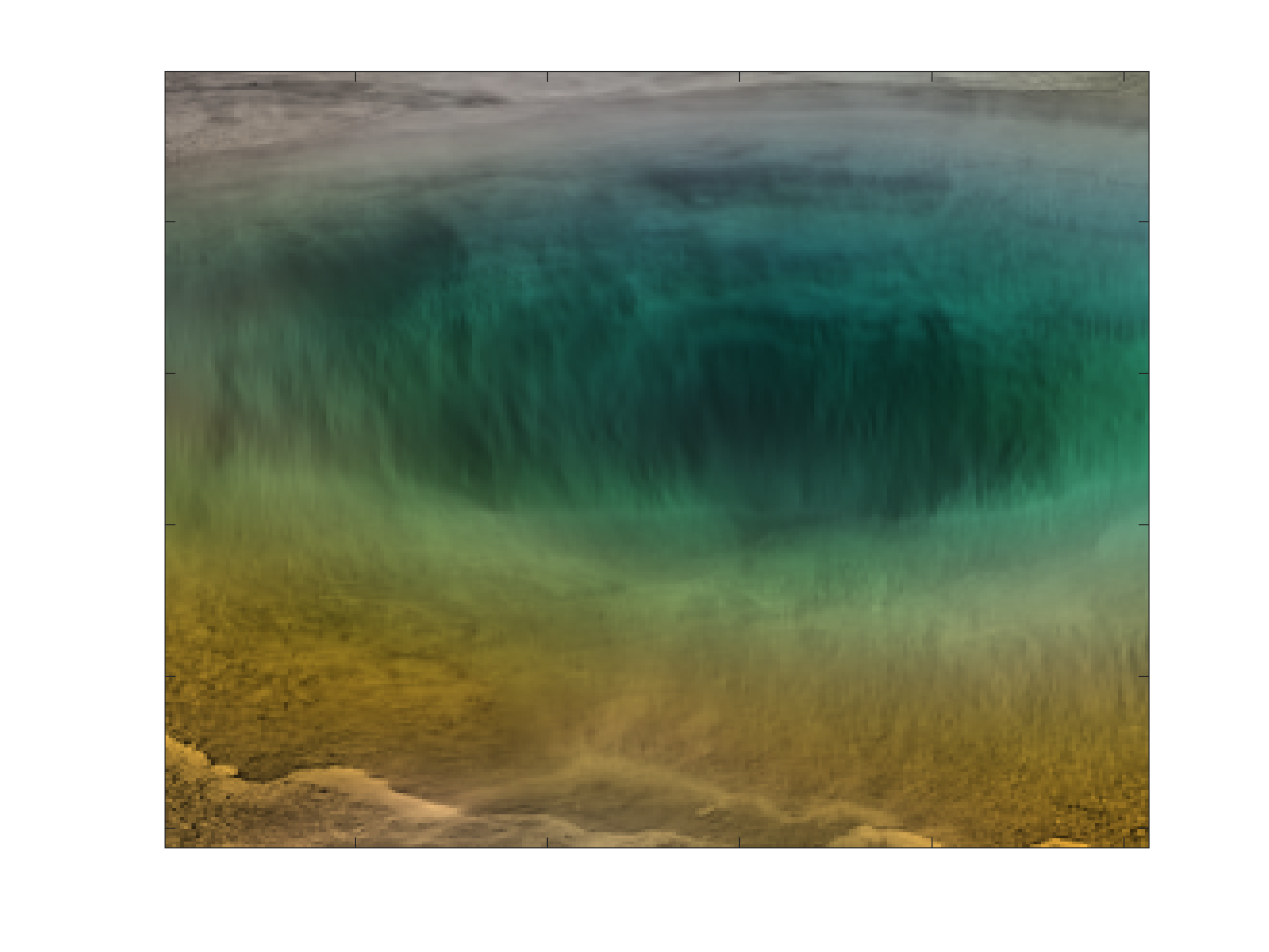}}
\subfloat[$t=0.6$]{\includegraphics[width=0.35\textwidth]{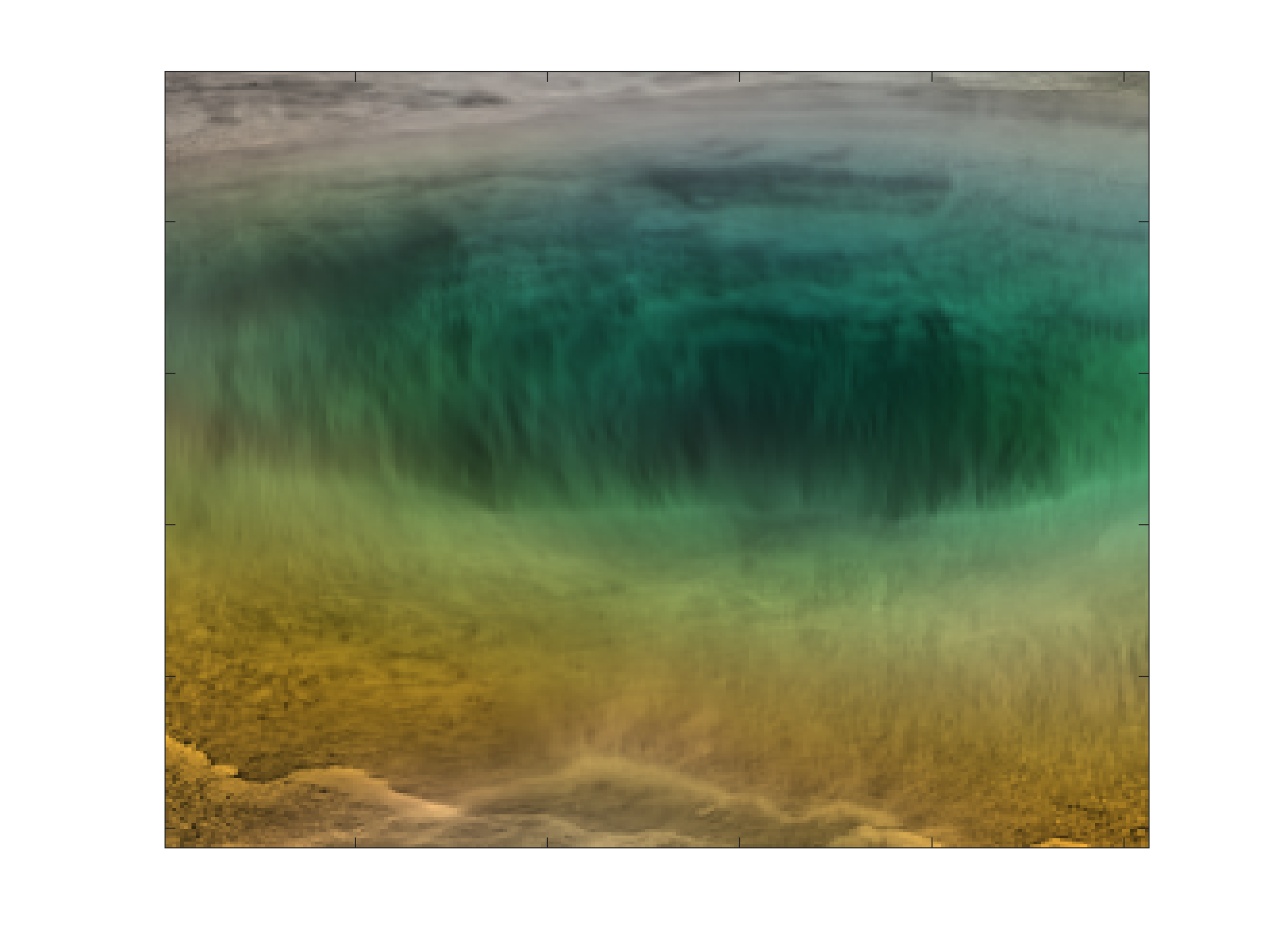}}\\
\subfloat[$t=0.7$]{\includegraphics[width=0.35\textwidth]{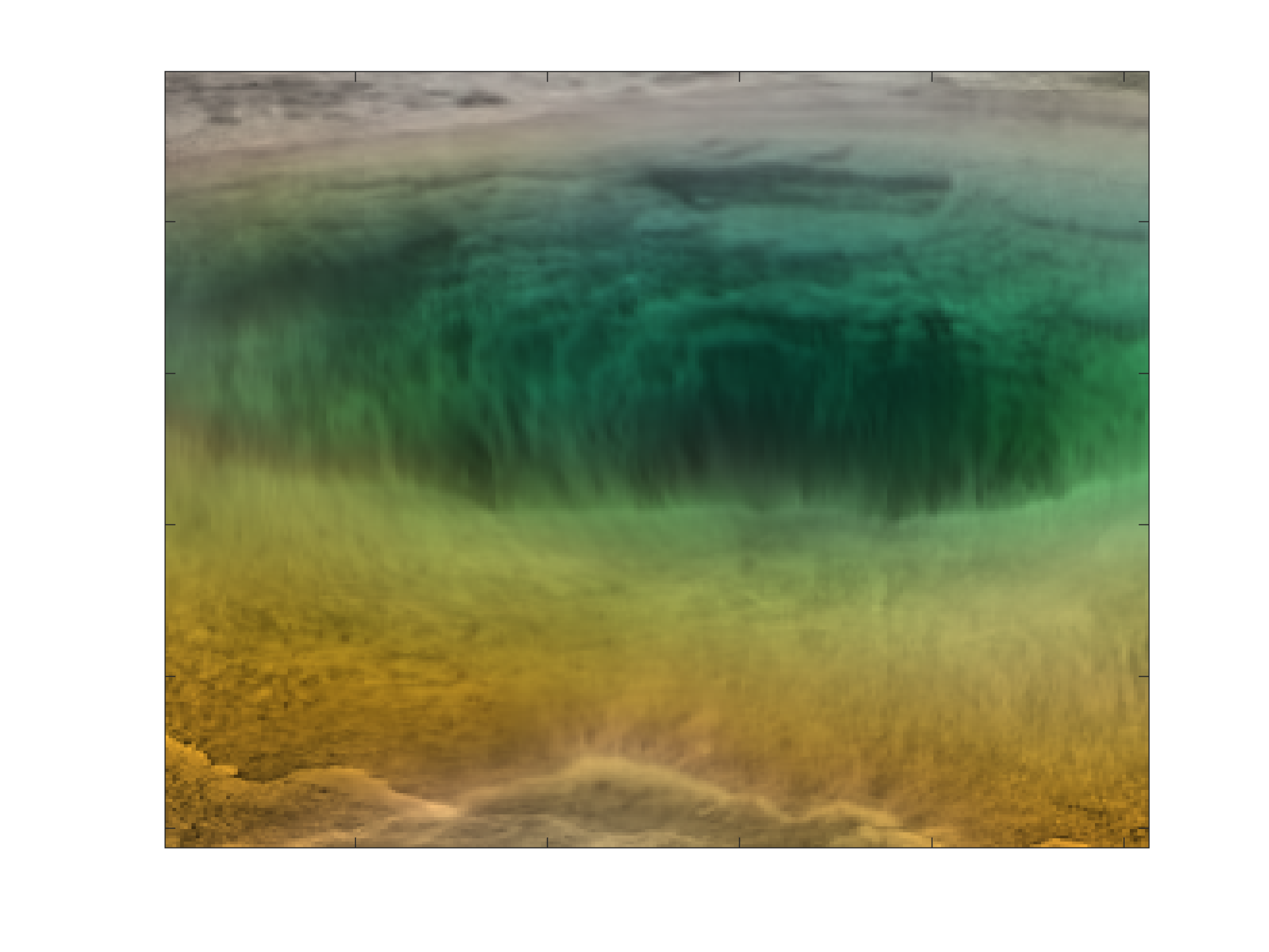}}
\subfloat[$t=0.8$]{\includegraphics[width=0.35\textwidth]{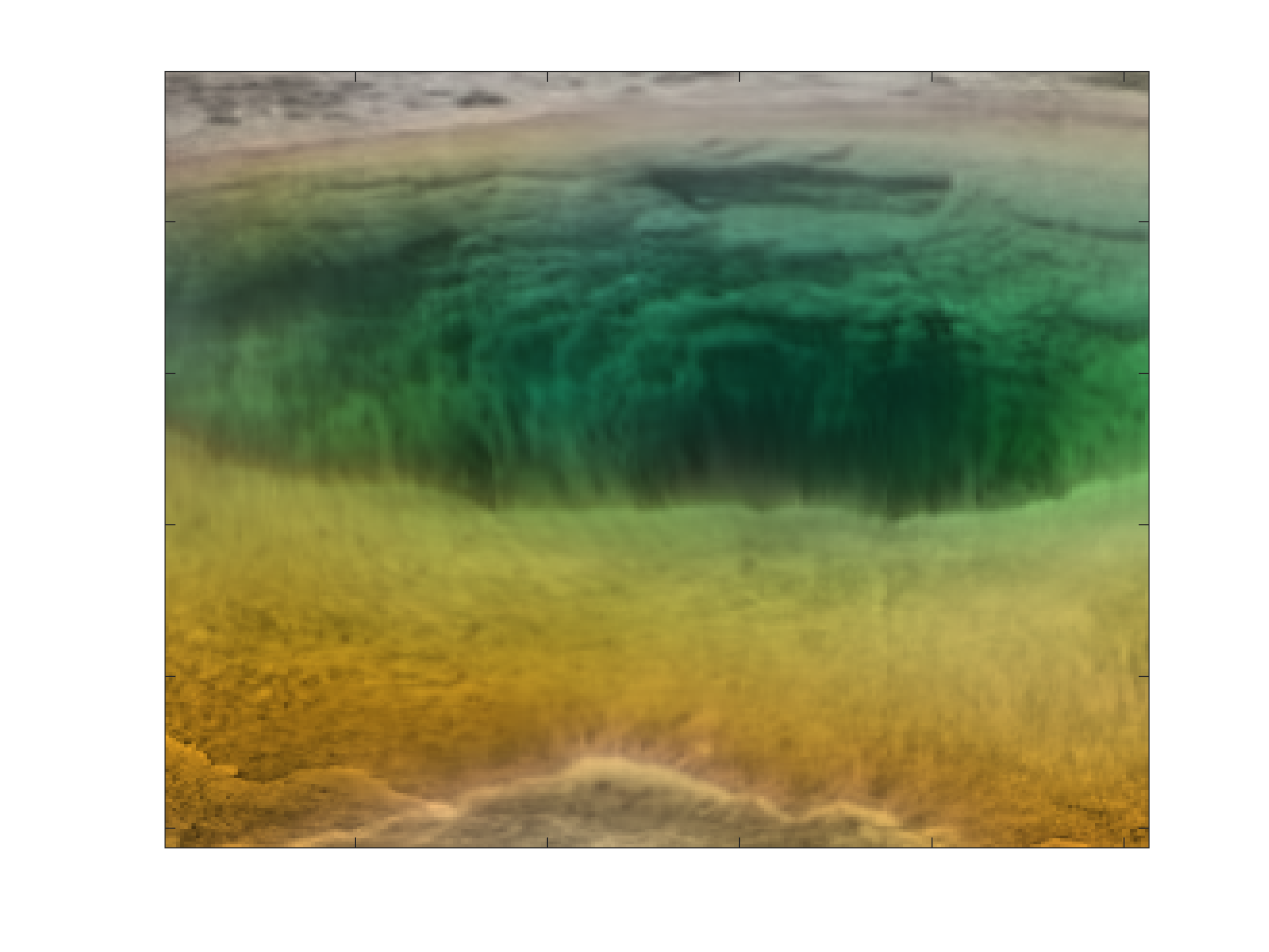}}
\subfloat[$t=0.9$]{\includegraphics[width=0.35\textwidth]{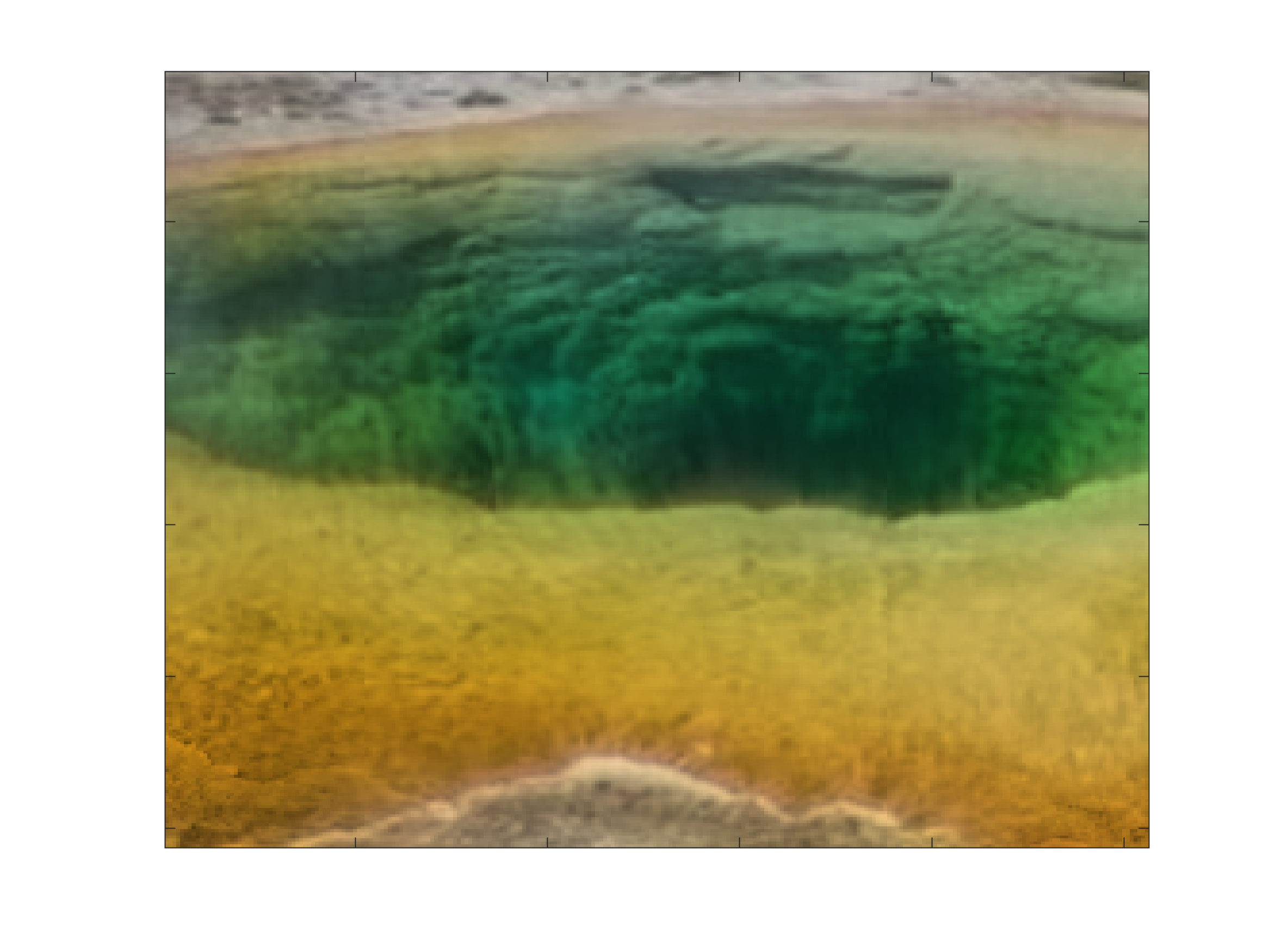}}
\caption{Interpolation with $\gamma=0.3$}
\label{fig:ex3interp}
\end{figure}

As alluded to previously, the different components of a image may stand for different color channels. For instance, a color image can be viewed as a vector-valued density with three components that represent red (R), green (G), blue (B), respectively. Thus, it is straightforward to use vector-valued optimal mass transport to compare and interpolate such color images. Below we explain three representative examples shown in Figures~\ref{fig:ex1marginals} through \ref{fig:ex3interp}, that highlight the mechanism of vector-valued transport.

First consider the two color images ($64 \times 64$) shown in Figure~\ref{fig:ex1marginals}. The intensity in each is a Gaussian distribution centered at a different location. The two distributions are of different color, thereby the corresponding vectorial-valued mass is distributed differently across the three components. The result of interpolating between the two, with $\gamma=0.001$, is shown in Figure~\ref{fig:ex1interp}. As we can see from the subplots, the displacement of the mass appears to run at constant speed between $\rho_0$ to $\rho_1$ while, at the same time, the color is changing gradually as mass flows between the vectorial components.

Figure~\ref{fig:ex2marginals} shows yet another example of a similar nature. The initial density is centered an it is white, which signifies equally mass distribution across the three color channels/components. The terminal density on the other hand has four separated masses of different color. The dimensions of the images are $128$ by $128$. The density flow shown in Figure~\ref{fig:ex2interp}, based on our technique with $\gamma=0.01$, smoothly interpolates by dispacing the intensity and color profiles in a seemingly natural manner.

Finally, in Figures~\ref{fig:ex3marginals} and \ref{fig:ex3interp} we display the result of interpolating real-life images. The marginal distributions shown in Figure~\ref{fig:ex3marginals} are two photos ($256 \times 256$) of two geothermal basins in Yellowstone Park, where bacterial growth give them distinctly different colors and hues. The result of interpolating the corresponding vector-valued distributions is depicted in Figure~\ref{fig:ex3interp}, taking $\gamma=0.3$. The flow of images produces a sequence of natural looking images transitioning from one to the other.

In all the examples, we observe the apparently natural displacement
of intensity and color that should be contrasted with potentially undesirable ``push-pop'' effects of linear interpolation.

\section{Conclusions and further research} \label{sec:conclusion}

Our early motivation, as noted in the introduction, has been to devise a suitable geometry to study flows of probability or power distribution in problems of signal analysis and fusion of vectorial data. However, the framework, very much as the broader subject of optimal mass transport, has application in a wider range of ideas. In particular, the connection between transport geometry and properties of an underlying space (e.g., curvature in the Bakry-Emery theory) may have important implications here as well. More specifically, we are interested in applying this methodology to studying the robustness of various networks as was done in \cite{Sandhu, Sandhu1} for biological and financial networks, and \cite{WanJonBan13} for communications networks.

\subsubsection{Biological networks}
The study of cellular networks (e.g., signalling and transcription) has become a major enterprise in systems biology; see \cite{Alo06} and the references therein.
One of the key problems is understanding global properties of cellular networks, in order to differentiate a diseased state from a normal cellular state. As is argued in several places \cite{DemMan05,WesBiaSevTes12,Sandhu} network properties may help in formulating systems biological concepts that could lead to novel therapies for a number of diseases including cancer. This would involve integrating genetic, epigenetic, and protein-protein interaction networks.

\subsubsection{Financial networks}
Stock-data and financial transactions provide an insight into the vast globe-wide financial network of human activities. The health of the national and world economy is reflected in the robustness and self-regulatory properties of the markets. Long range correlations are responsible for cascade failures due to financial insolvency. Indeed, multiple exposures of companies is often the root cause of infectious propagation of balance sheet insolvency with catastrophic effects. It is of interest to understand the relation between the various financial parameters (assets, liability, capital) that quantify the stress and the buffer capabilities of financial institutions with the network connectivity and interdependence (weighted network Laplacian) so as to assess risk of cascade failures, fragility, and devise ways to mitigate such effects. See \cite{Sandhu1} and the references therein.

At closer inspection, many of the aforementioned problem areas involve finer attributes of the studied objects, which may be more suitably treated and studied as vector-valued distributions. Thus, we hope that the present work provides a starting point for such an endeavor.

\section{Acknowledgements}
This project was supported by AFOSR grants (FA9550-15-1-0045 and FA9550-17-1-0435), grants from the National Center for Research Resources (P41-
RR-013218) and the National Institute of Biomedical Imaging and Bioengineering (P41-EB-015902), National Science Foundation (NSF), and grants from National Institutes of Health (P30-CA-008748, 1U24CA18092401A1, R01-AG048769).

\newpage
\bibliographystyle{IEEEtran}
\bibliography{refs}

\end{document}